\documentclass[times,final]{elsarticle}

\usepackage{amsmath} 
\usepackage{amsthm} 

\usepackage{jcomp}
\usepackage{framed,multirow}

\usepackage{amssymb}
\usepackage{latexsym}

\usepackage{url}
\usepackage[dvipsnames]{xcolor}

\usepackage{subcaption} 
\usepackage{hyperref}   
\usepackage{siunitx}    

\journal{Journal of Computational Physics}

\newcommand\rone[1]{\textcolor{black}{#1}}
\newcommand\rtwo[1]{\textcolor{black}{#1}}

\definecolor{mycustomred}{HTML}{eb3b5a}
\definecolor{mycustomgreen}{HTML}{20bf6b}
\newcommand\red[1]{\textcolor{mycustomred}{#1}}
\newcommand\green[1]{\textcolor{mycustomgreen}{#1}}

\begin{document}

\verso{Given-name Surname \textit{etal}}

\begin{frontmatter}

\title{A robust hybridizable discontinuous Galerkin scheme with harmonic averaging technique for steady state of real-world semiconductor devices}%

\author[1,2]{Qingyuan \snm{Shi}}
\cortext[cor1]{Corresponding author: 
  chijie@tsinghua.edu.cn (Chijie Zhuang);  
  linbo@u.nus.edu (Bo Lin)}
\author[3]{Yongyong \snm{Cai}}
\author[1,2]{Chijie \snm{Zhuang}\corref{cor1}}
\author[4]{Bo \snm{Lin}\corref{cor1}}
\author[2]{Dan \snm{Wu}}
\author[1]{Rong \snm{Zeng}}
\author[4]{Weizhu \snm{Bao}}

\address[1]{Department of Electrical Engineering, Tsinghua University, Beijing 100084, China}
\address[2]{Beijing Huairou Laboratory, Beijing 101400, China}
\address[3]{Laboratory of Mathematics and Complex Systems and School of Mathematical Sciences, Beijing
Normal University, Beijing 100875, China}
\address[4]{Department of Mathematics, National University of Singapore, Singapore 119076, Singapore}

\received{}
\finalform{}
\accepted{}
\availableonline{}
\communicated{}

\begin{abstract}
Solving real-world nonlinear semiconductor device problems modeled by the drift-diffusion equations coupled with the Poisson equation (also known as the Poisson-Nernst-Planck equations) necessitates an accurate and efficient numerical scheme which can avoid non-physical oscillations even for problems with extremely sharp doping profiles. In this paper, we propose a flexible and high-order hybridizable discontinuous Galerkin (HDG) scheme with harmonic averaging (HA) technique to tackle these challenges. The proposed HDG-HA scheme combines the robustness of finite volume Scharfetter-Gummel (FVSG) method with the high-order accuracy and $hp$-flexibility offered by the locally conservative HDG scheme. The coupled Poisson equation and two drift-diffusion equations are simultaneously solved by the Newton method. Indicators based on the gradient of net doping $N$ and solution variables are proposed to switch between cells with HA technique and high-order conventional HDG cells, utilizing the flexibility of HDG scheme.
Numerical results suggest that the proposed scheme does not exhibit oscillations or convergence issues, even when applied to heavily doped and sharp PN-junctions. Devices with circular junctions and realistic doping profiles are simulated in two dimensions, qualifying this scheme for practical simulation of real-world semiconductor devices.

\end{abstract}


\end{frontmatter}


\section{Introduction}
\begin{figure}[!ht]
\centering
\centering
\includegraphics[width=0.72\textwidth]{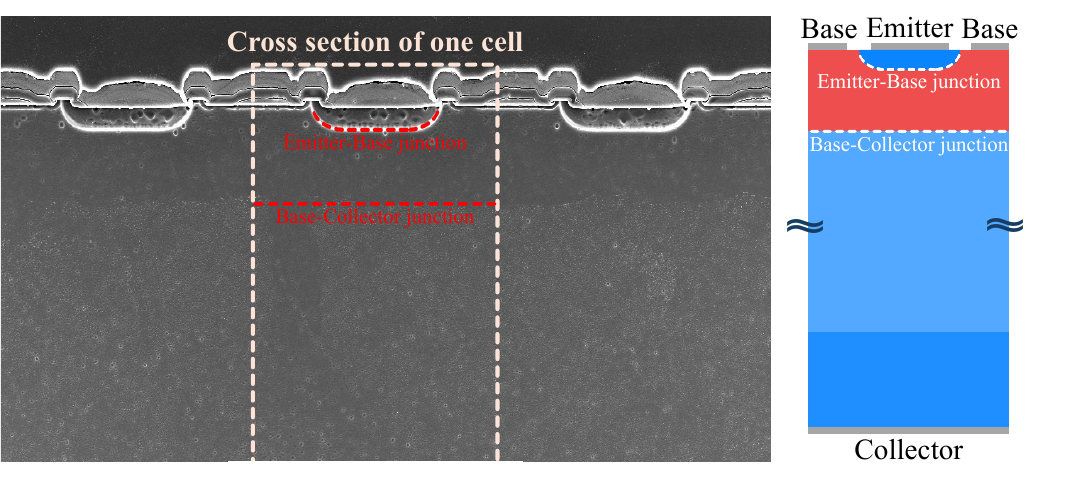}
\caption{\rone{Structure of a power BJT: SEM cross-section of multiple cells (left) and schematic representation of a single cell (right).}}
\label{fig:real_BJT_illustration}
\end{figure}

The simulation of semiconductor devices requires solving for the interaction between electric field and charge transport, modeled by the Poisson-Nernst-Planck (PNP) equations \cite{Selberherr1984Analysis}. Although numerical methods on PNP equations have been studied for decades, the essential difficulties of semiconductor problems, where the two types of equations are strongly coupled, still challenge the robust and accurate numerical algorithms for semiconductor devices, especially when the power rating and voltage level of modern power semiconductor devices have substantially increased in recent years \cite{Zeng2019IGCT}. The major obstacles for an ideal numerical scheme are: 

\begin{enumerate}
    \item Strongly nonlinear coupling: The strong coupling between Poisson and convection-diffusion equations in real-world devices' operating conditions necessitates solving for a coupled nonlinear system of equations instead of Gummel-like decoupled and linearized iteration techniques \cite{Mock1972GummelConvergence}.
    
    \item Numerical stability near junctions: Across the junction structures in semiconductor devices, net doping concentration present sharp jumps within thin junction layers, resulting in a depletion region with high electric field and rapid gradient of carrier concentrations over tens of orders of magnitude. The convection-dominated nature of the depletion region poses numerical challenge for stable and robust numerical schemes.
    
    \item Local conservation: The accurate reconstruction of electric field and current density inside a device, with which the avalanche generation rate and overheating point can be computed, plays a vital role in practical device analysis and optimization. Therefore the numerical scheme should be locally conservative.

    \item Characteristics of realistic power devices: Real-world power semiconductor devices, with high doping concentrations (more than $10^{20} {\rm cm}^{-3}$), narrow junction layers and high voltage blocking standard, demand more robust and efficient algorithms. \rone{The structure of a power bipolar junction transistor (BJT) observed through a scanning electron microscope (SEM), is illustrated in Fig. \ref{fig:real_BJT_illustration}.} Many of the conventional algorithms fail to work for huge doping jumps. The dielectric relaxation time restriction \cite{Lin2020JCP} under high doping would prohibit conventional techniques of solving steady state of devices using explicit time marching \cite{Liu2004LDG1, Liu2007LDG2}.
\end{enumerate}

In the past decades, many numerical schemes have been proposed to address the above issues, here we categorize the majorities of them into three types: finite volume method (FVM), finite element method (FEM) and discontinuous Galerkin (DG) scheme.

The finite volume method accompanied with the Scharfetter-Gummel (SG) discretization for drift-diffusion (DD) equations \citep{ScharfetterGummel1969}, known as the FVSG scheme \cite{Bank1983FVSG}, has become the industrial standard for semiconductor simulation and widely used in commercial simulators. It is both robust and flexible on unstructured meshes, guaranteeing local conservation properties while capturing sharp layers at junctions without non-physical oscillations. Albeit its advantages, the inherited low order nature of the FVSG scheme makes it inaccurate, requiring fine meshes which can be too computationally burdensome for three dimensional simulations. Although high-order SG schemes have been proposed \cite{Nguyen2022SGCC}, its construction is complicated and it is hard to be generalized to unstructured meshes. Other types of FVM like element edge based FVM \cite{Sanchez2021EEB_SG} and cell-centered FVM \cite{Rupp2014CellCenterFVSG} have been proposed, but the same SG discretization was adopted for fluxes, resulting in the same low accuracy.

The conventional finite element methods (FEM) can achieve higher-order accuracy than the FVSG scheme for general problems on high-order polynomial spaces. However, these methods suffer from non-physical oscillations and non-convergent issues when applied to the drift-diffusion equations of semiconductors \cite{SongWang1999Survey}. These oscillations and convergence difficulties primarily arise from the nonlinear and convection-dominated nature of semiconductor problems. Therefore, special stabilization techniques have to be imposed to get convergent results. The streamline upwind/Petrov-Galerkin (SUPG) method \cite{Brooks1982SUPG}, which introduces stabilization terms to enhance stability and reduce oscillations, has been extensively studied on convection-diffusion and semiconductor equations \cite{Swaminathan1992SUPG_CVFEM, Di2018SUPG, SandiaLin2010SUPGParallel, Wang2021_ZJU_SUFVFEM}. However, there are no general choices of their stabilization parameters, and they cannot completely eliminate oscillations either \cite{Bochev2013CVFEMSG}. Apart from adding artificial stabilization terms, numerous FEM schemes with harmonic averaging or exponential fitting techniques have been studied. Edge-averaged finite element (EAFE) \cite{Xu1999MonotoneEAFE, Sacco2015EAFE3D}, mixed finite element method \cite{Brezzi1989ExpFit, Miller1994TetrahedralMixFEM} and inverse average type finite element methods \cite{Markowich1988InverseAverage, Zhang2022IAFEM, Lazarov2012ExpFitFEM}, which all utilize harmonic averaging techniques to produce M-matrices, are shown to be stable. Changing basis functions from polynomials to exponential functions \cite{Wang1999ExpBasis, Wang1997ExpBasis} also showed its stability on convection-dominated problems. However, all the aforementioned FEM schemes lose the local conservation property, preventing them from being applied to real-world simulation and analysis of devices. The control volume finite element method with multi-dimensional Scharfetter–Gummel upwinding (CVFEM-SG) \cite{Bochev2013CVFEMSG, Li2022_ZJU_CVFEM} method, which utilizes the SG discretization of current fluxes, offers the advantages of stability and local conservation. However, like the FVSG scheme, it still suffers from a low order of accuracy. 

Combining the advantages of the FVM and FEM schemes, the discontinuous Galerkin schemes \cite{CockburnShu2012DG_Book} realize local conservation and high-order accuracy on general unstructured meshes, showing potential advantages for semiconductor devices simulation. A large number of studies have applied DG methods to solve convection-dominated problems \cite{Cockburn1998ConvectionDominatedDG}, singularly perturbed problems and PNP problems. Specifically, DG-based methods such as the Runge-Kutta discontinuous Galerkin (RKDG) method \cite{Shu1995MixedRKDGSemiconductor}, the local discontinuous Galerkin (LDG) method \cite{Liu2004LDG1, Liu2007LDG2, Liu2016AnalysisLDG} and the hybridizable discontinuous Galerkin (HDG) method \cite{Chen2019HDG, Feng2023HDG_ZJU} have been employed for semiconductor simulations. However, most of these studies used only explicit temporal discretization \cite{Shu1995MixedRKDGSemiconductor, Chen2020AccessDGSemiconductpr} and decoupled iteration method \cite{Chen2020AccessDGSemiconductpr, Feng2023HDG_ZJU} for Poisson and drift-diffusion equations, making it easier for implementations. Moreover, most existing studies only simulated simple cases where no PN-junction \cite{Chen2019HDG} or only moderate junctions (where doping only slightly changes across the junction, like high-low junctions in pure N-type semiconductors \cite{Liu2004LDG1, Liu2007LDG2}) present, thereby substantially reducing computational difficulty from real-world problems. To the best of our knowledge, few studies have conducted on realistic huge and extremely sharp doping profiles using an implicit simulation of coupled nonlinear semiconductor equations with DG methods, yet no studies have ever revealed their convergence properties of problems from the industry.

In this paper, we focus on the practical simulation of real-world power semiconductor devices with the following contributions. We adapt the pivotal harmonic-averaging technique of the FVSG scheme to the flexible HDG scheme, forming the hybridizable discontinuous Galerkin scheme with harmonic average technique (HDG-HA scheme) that combines the merits of both methods. Furthermore, we propose straightforward but effective indicators to adaptively switch between harmonic averaging cells and conventional HDG cells and show primary results on the $hp$-adaptability of the HDG-HA scheme. Finally, we apply the HDG-HA scheme on real-world power semiconductor devices, from which the blocking state and punch-through effects of devices are analyzed withstanding thousands of volts.

The rest of this paper is organized as follows. Section \ref{sec:ModelsAndDiscretization} provides a brief description of the physical and mathematical model of semiconductor devices. In Section \ref{sec:HDG}, we present the hybridizable discontinuous Galerkin for semiconductor devices' PNP equations. In Section \ref{sec:HA}, we elaborate the construction of the proposed harmonic averaging technique on unstructured meshes. In Section \ref{sec:Experiments}, we conduct extensive experiments to validate the convergence and accuracy of the proposed HDG-HA scheme, as well as the indicators for the HA-adaptivity. The possible future extension of $hp$-adaptivity with HA-adaptivity is also illustrated. In Section \ref{sec:RealWorldExamples}, two real-world power semiconductor devices are simulated, analyzing the current-voltage characteristics of a PiN power diode and the punch-through effect of a bipolar junction transistor, respectively.

\section{The model equations and their nondimensionalization}
\label{sec:ModelsAndDiscretization}

\subsection{The governing equations of semiconductor devices in two dimensions (2D)}

Electrostatic potential $\psi:=\psi(\mathbf{x},t)\in\mathbb{R}$ and carrier densities $n:=n(\mathbf{x},t), p=p(\mathbf{x},t)\in\mathbb{R}$  ($n$ for electrons and $p$ for holes) in semiconductor material domain $\Omega \subset \mathbb{R}^2$ are modeled by the Poisson equation coupled with drift-diffusion equations (i.e., the Poisson-Nernst-Planck equations)
\begin{equation}
    \label{GoverningEquations}
    \begin{cases}
        \nabla^2 \psi(\mathbf{x},t)+\frac{q}{\varepsilon}\left(p(\mathbf{x},t)-n(\mathbf{x},t)+N(\mathbf{x})\right)=0, \\
        \partial_t n(\mathbf{x},t)- \frac{1}{q}\nabla \cdot \mathbf{J}_{n}(\mathbf{x},t) + R_{n}(n,p)=0, \\
        \partial_t p(\mathbf{x},t)+\frac{1}{q}\nabla \cdot \mathbf{J}_{p}(\mathbf{x},t) + R_{p}(n,p)=0, \\
    \mathbf{J}_n(\mathbf{x},t)=-q\,\mu_n\,n(\mathbf{x},t)\,\nabla \psi(\mathbf{x},t) + q\,D_n \nabla n(\mathbf{x},t), \\
         \mathbf{J}_p(\mathbf{x},t)=-q\,\mu_p\,p(\mathbf{x},t)\,\nabla \psi(\mathbf{x},t) - q\,D_p \nabla p(\mathbf{x},t),
    \end{cases} \mathbf{x} \in \Omega,\quad t>0.
\end{equation}
Here $\nabla^2=\Delta$ is the Laplace operator, $\nabla$ is the gradient operator, $\nabla\cdot$ is the divergence operator, $\mathbf{J}_n\in\mathbb{R}^2$ and $\mathbf{J}_p\in\mathbb{R}^2$ are the currents,  $R_n\in\mathbb{R}$ and $R_p\in\mathbb{R}$ are the source terms of electrons and holes, respectively, and $N:=N(\mathbf{x})\in\mathbb{R}$ is the net doping concentration, which is pre-defined by the manufacturing process.  Inside $\Omega$, semiconductor dielectric constant $\varepsilon$, Boltzmann constant $k_{\rm b}$, elementary charge $q$, device temperature $T$, carrier mobilities $\mu_n, \mu_p$ and carrier diffusion coefficients $D_n,D_p$ are all considered as constants. In addition, we assume the Einstein relationship $D_n=\mu_n V_T,D_p=\mu_p V_T$ holds, where the thermal voltage is defined by  $V_T = k_{\rm b}T/q$.   

 For the source terms $R_n,R_p$, Shockley-Read-Hall (SRH) and Auger recombinations \citep{Selberherr1984Analysis} are included as
\begin{equation}
    \label{RecombinationTerms}
    \begin{aligned}
        & R_n(n,p) =R_p(n,p)=R_{\rm SRH}(n,p)+R_{\rm Auger}(n,p),  \\
        & R_{\rm SRH}(n,p) = \frac{n p-n_{\rm ie}^2}{\tau_p\left(n+n_{\rm ie}\right)+\tau_n\left(p+n_{\rm ie}\right)},\quad R_{\rm Auger}(n,p) = \left(C_n n+C_p p\right)\left(n p-n_{\rm ie}^2\right),
    \end{aligned}
\end{equation}
where $n_{\rm ie}$ is the effective intrinsic carrier density in the semiconductor material, $C_n$ and $C_p$ are the corresponding recombination coefficients, $\tau_n$ and $\tau_p$ are the electron and hole lifetimes, respectively. In this paper, all the physical constants and parameters take the values of silicon material at $\SI{300}{K}$ temperature, and we refer to the Appendix for details.

The boundary of the semiconductor domain $\partial \Omega$ is composed of two parts: the Dirichlet boundary, $\partial \Omega_{\rm Dirichlet}$, and the homogenous Neumann boundary, $\partial \Omega_{\rm Neumann}$. The Dirichlet boundary corresponds to interfaces where the semiconductor material forms Ohmic contacts with metal electrodes. And the homogenous Neumann boundary is characterized by the semiconductor material adjoining insulator boundaries, where there are no surface charges or current fluxes. For Ohmic contacts (i.e., Dirichlet boundary) under Maxwell-Boltzmann statistics, boundary conditions are given by 
\begin{equation}
    \label{DirichletBoundary}
    \begin{aligned}
        & \psi|_{\partial \Omega_{\rm Dirichlet}}= V_{\rm applied} + V_T \ln\left(\frac{n}{n_{\rm ie}}\right),
        \quad n|_{\partial \Omega_{\rm Dirichlet}} = \frac{1}{2}\left(N+\sqrt{N^2+4n_{\rm ie}^2}\right),
        \quad p|_{\partial \Omega_{\rm Dirichlet}} = \frac{1}{2}\left(-N+\sqrt{N^2+4n_{\rm ie}^2}\right),
    \end{aligned}
\end{equation}
where $V_{\rm applied}$ represents the applied external voltage, which may vary across different electrodes. For the Neumann boundaries, no-flux conditions are assumed
\begin{equation}
    \label{NeumannBoundary}
    \begin{aligned}
        & \nabla\psi|_{\partial \Omega_{\rm Neumann}} \cdot \mathbf{n}=0, \quad \mathbf{J}_n|_{\partial \Omega_{\rm Neumann}} \cdot \mathbf{n} =0, \quad \mathbf{J}_p|_{\partial \Omega_{\rm Neumann}} \cdot \mathbf{n} =0,
    \end{aligned}
\end{equation}
where $\mathbf{n}$ is the unit outer normal vector of $\partial \Omega_{\rm Neumann}$.
\subsection{Nondimensionless form}
\label{subsec:Nondimensionless}

We first scale the governing equations \eqref{GoverningEquations} using the singular perturbation scaling \cite{Peter2013StationaryEquations} to get the dimensionless equations, i.e. $\psi=\tilde{\psi}V^*$, $N=\tilde{N}N^*$, $n=\tilde{n}N^*$, $p=\tilde{p}N^*$,  $\mathbf{J}_n=\tilde{\mathbf{J}}_nJ^*$, $\mathbf{J}_p=\tilde{\mathbf{J}}_pJ^*$, $\mu_n=\tilde{\mu}_n\mu^*$, $\mu_p=\tilde{\mu}_p\mu^*$, $D_n=\tilde{D}_nD^*$, $D_p=\tilde{D}_pD^*$, $\mathbf{x}=\tilde{\mathbf{x}}x^*$ ($\Omega\to \tilde{\Omega} x^*$), $t=\tilde{t}t^*$, for the Dirichlet boundary conditions, $V_{\rm applied}=\widetilde{V}_{\rm applied} V^*$, and for the source terms  $R_n=R^*\tilde{R}_n$, $R_p=R^*\tilde{R}_p$, $n_{\rm ie}=N^*\tilde{n}_{\rm ie}$, $\tau_p=t^*\tilde{\tau}_p$, $\tau_n=t^*\tilde{\tau}_n$, $C_n=C^*\tilde{C}_n$, $C_p=C^*\tilde{C}_p$,  where $V^*$, $N^*$, $J^*$, $\mu^*$, $D^*$, $x^*$, $t^*$, $R^*$ and $C^*$ are the corresponding base units. 
\begin{itemize}
    \item Potential $\psi$ is scaled by thermal voltage $V^*=V_T$.
    \item Carrier concentrations $n,p$ are scaled by the maximum doping concentration $N^*=\max(N)$.
    \item Length scale $x^*$ is typically chosen as the same order to the diameter of $\Omega$, i.e. $x^*=\text{diameter}(\Omega)$. 
    \item Carrier mobility $\mu_n,\mu_p$  and diffusion coefficients $D_n,D_p$ are scaled by their maximum value $D^*\approx \max(D_n,D_p)$, $\mu^*=\frac{D^*}{V^*}$, so that $\frac{\mu_n}{\mu^*}\sim 1,\frac{\mu_p}{\mu^*}\sim 1,\frac{D_n}{D^*}\sim 1,\frac{D_p}{D^*}\sim 1$. 
    \item The remaining scaling base units depend on  $x^*,V^*,N^*,D^*$, i.e., $J^*=q\frac{D^*N^*}{x^*}, R^*=\frac{D^*N^*}{\left(x^*\right)^2}, t^*=\frac{\left(x^*\right)^2}{D^*}$, $C^*=\frac{1}{(N^*)^2t^*}$.
\end{itemize}
We list the typical scaling  units for silicon material in the Appendix. Plugging the scaled variables into \eqref{GoverningEquations} and removing~~$\tilde{}$, we obtain the dimensionless PNP equations
\begin{equation}
    \label{dimless}
    \begin{cases}
        \lambda^2\nabla^2 \psi(\mathbf{x},t)+\left(p(\mathbf{x},t)-n(\mathbf{x},t)+N(\mathbf{x})\right)=0, \\
        \partial_t n(\mathbf{x},t)-\nabla \cdot \mathbf{J}_{n}(\mathbf{x},t) + R_{n}(n,p)=0, \\
        \partial_t p(\mathbf{x},t)+\nabla \cdot \mathbf{J}_{p}(\mathbf{x},t) + R_{p}(n,p)=0, \\
        \mathbf{J}_n(\mathbf{x},t)=-\mu_n\,n(\mathbf{x},t)\,\nabla \psi(\mathbf{x},t) + D_n \nabla n(\mathbf{x},t), \\
        \mathbf{J}_p(\mathbf{x},t)=-\mu_p\,p(\mathbf{x},t)\,\nabla \psi(\mathbf{x},t) - D_p \nabla p(\mathbf{x},t),    
    \end{cases} \mathbf{x} \in \Omega,\quad t>0, 
\end{equation}
 where the parameter $\lambda>0$ in the Poisson equation is given by $\lambda=\sqrt{\frac{\varepsilon V^*}{qN^*\left(x^*\right)^2}}$. The Neumann boundary conditions are  in the same form as \eqref{NeumannBoundary}, and the Dirichlet boundary conditions become
 \begin{equation}
    \label{DirichletBoundaryRes}
    \begin{aligned}
        & \psi|_{\partial \Omega_{\rm Dirichlet}}= V_{\rm applied} +  \ln\left(\frac{n}{n_{\rm ie}}\right), \quad n|_{\partial \Omega_{\rm Dirichlet}} = \frac{1}{2}\left(N+\sqrt{N^2+4n_{\rm ie}^2}\right),  \quad p|_{\partial \Omega_{\rm Dirichlet}} = \frac{1}{2}\left(-N+\sqrt{N^2+4n_{\rm ie}^2}\right).
    \end{aligned}
\end{equation}
 In addition, the dimensionless recombination source terms are defined as
 \begin{equation}
\label{ScaledRecombinationRate}
\begin{aligned}
 & R_n(n,p) =R_p(n,p)=R_{\rm SRH}(n,p)+R_{\rm Auger}(n,p), \\
 & R_{\rm SRH}(n,p) = \frac{n p-n_{\rm ie}^2}{\tau_p\left(n+n_{\rm ie}\right)+\tau_n\left(p+n_{\rm ie}\right)},\quad R_{\rm Auger}(n,p) = \left(C_n n+C_p p\right)\left(n p-n_{\rm ie}^2\right),
\end{aligned}
 \end{equation}
where $\tau_n, \tau_p$ are the dimensionless lifetime of carriers, $C_n$ and $C_p$ are the dimensionless recombination coefficients for electrons and holes, respectively, and $n_{\rm ie}$ is the dimensionless effective intrinsic concentration.
 
 In one dimension (1D), we assume the system is homogeneous along the rest direction, and the semiconductor device equations \eqref{dimless} remain the same for $\Omega\subset\mathbb{R}$. In this paper, we shall focus on the steady states of \eqref{dimless} for $\Omega\subset\mathbb{R}^d$ ($d=1,2$), i.e. the time-independent PNP equations for $p(\mathbf{x}),n(\mathbf{x}),\psi(\mathbf{x})\in\mathbb{R}$, and $\mathbf{J}_n(\mathbf{x}),\mathbf{J}_p(\mathbf{x})\in\mathbb{R}^d$,
\begin{equation}
    \label{steady}
    \begin{cases}
        \lambda^2\nabla^2 \psi(\mathbf{x})+\left(p(\mathbf{x})-n(\mathbf{x})+N(\mathbf{x})\right)=0, \\
        -\nabla \cdot \mathbf{J}_{n}(\mathbf{x}) + R_{n}(n,p)=0,\\
        \nabla \cdot \mathbf{J}_{p}(\mathbf{x}) + R_{p}(n,p)=0,\\
        \mathbf{J}_n(\mathbf{x})=-\mu_n\,n(\mathbf{x})\,\nabla \psi(\mathbf{x}) + D_n \nabla n(\mathbf{x}),\\
        \mathbf{J}_p(\mathbf{x})=-\mu_p\,p(\mathbf{x})\,\nabla \psi(\mathbf{x}) - D_p \nabla p(\mathbf{x}), \\
    \end{cases}
    \mathbf{x} \in \Omega,
\end{equation}
where $R_{n}(n,p)$ and $R_{p}(n,p)$ are given in \eqref{ScaledRecombinationRate} and the boundary conditions are prescribed in \eqref{NeumannBoundary} and \eqref{DirichletBoundaryRes}.
 
\section{The hybridizable discontinuous Galerkin scheme}
\label{sec:HDG}

In this section, we present the hybridizable discontinuous Galerkin scheme for solving the time independent semiconductor device equations \eqref{steady}.

\subsection{Discretization of the hybridizable discontinuous Galerkin scheme}

In $d$ dimensional space $\Omega\subset\mathbb{R}^d$, by introducing an auxiliary variable $\mathbf{E}=-\nabla\psi$, the time independent governing equations \eqref{steady} can be written as

\begin{equation}
    \label{FirstOrderGoverningEquations}
    \left\{
        \begin{aligned}
        &\lambda^2\nabla \cdot \mathbf{E}+n-p-N = 0, \\ 
        & - \nabla \cdot \mathbf{J}_{n}+R_n(n,p) = 0, \\ 
        &\nabla \cdot \mathbf{J}_{p}+R_p(n,p) = 0, \\    
        &\mathbf{E}+\boldsymbol{\nabla} \psi = 0, \\ 
        &\mathbf{J}_n-\mu_n n \mathbf{E}-D_n \boldsymbol{\nabla} n = 0,\\ 
        &\mathbf{J}_p-\mu_p p \mathbf{E}+D_p \boldsymbol{\nabla} p = 0.
        \end{aligned}
    \right.
\end{equation}

To introduce the HDG scheme, we adopt the same notations as \citep{CockburnNonlinear2009}. Let $\mathcal{T}_h$ represent the subdivision of $\Omega\subset\mathbb{R}^2$ such that either all elements $K \in \mathcal{T}_h$ are triangles or all $K \in \mathcal{T}_h$ are quadrilaterals. In 1D, $\mathcal{T}_h$ degenerates into a partition of an interval. The union of the edges of $K$ are denoted as $\partial \mathcal{T}_h:=\left\{ \partial K :K \in \mathcal{T}_h \right\}$. $\mathcal{S}^p(K)$ denotes the set of polynomials of total degree at most $p$ on $K$ if $K$ is a triangle,
or the set of polynomials of degree at most $p$ on $K$ if $K$ is a quadrilateral.
For any element $K \in \mathcal{T}_h$, we write $W^p(K)\equiv\mathcal{S}^p(K),\boldsymbol{V}^p(K)\equiv(\mathcal{S}^p(K))^d$. Now, we introduce the discontinuous ﬁnite element spaces

\begin{equation}
    \label{DGSpaces}
    \begin{aligned}
    & W_h^p=\left\{w \in L^2(\Omega):\left.w\right|_K \in W^p(K), \forall  K \in \mathcal{T}_h\right\}, \qquad \boldsymbol{V}_h^p=\left\{\boldsymbol{v} \in\left(L^2(\Omega)\right)^d:\left.\boldsymbol{v}\right|_K \in \boldsymbol{V}^p(K), \forall K \in \mathcal{T}_h\right\},
    \end{aligned}
\end{equation}
where $L^2(\Omega)$ is the standard space of square integrable functions on $\Omega$.

Let $\mathcal{E}^o_h$ and $\mathcal{E}^\partial_h$ be the set of interior and boundary element edges, respectively, with $\mathcal{E}_h$ being the union of $\mathcal{E}^o_h$ and $\mathcal{E}^\partial_h$. For $d=2$, the trace finite element space is defined as

\begin{equation}
    \label{TraceSpace}
    M_h^p=\left\{\mu \in L^2\left(\mathcal{E}_h\right):\left.\mu\right|_e \in \mathcal{S}^p(e), \forall e \in \mathcal{E}_h\right\},
\end{equation}
with $L^2\left(\mathcal{E}_h\right)$ the standard space of square integrable functions on $\mathcal{E}_h$. For $d=1$, the trace finite element space reduces to a vector space equipped with the discrete $l^2$ norm on discrete grid points.

For $K \in \mathcal{T}_h$, let $\left(\boldsymbol{w},\boldsymbol{v}\right)_K=\int_K \boldsymbol{w}\cdot \boldsymbol{v} \mathrm{d} \mathbf{x}$ be the inner product of two functions $\boldsymbol{w}, \boldsymbol{v} \in \left(L^2(K)\right)^d$, and $\left(u,v\right)_K=\int_K uv \mathrm{d}\mathbf{x}$ for $u, v \in L^2(K)$. Similarly for $e \in \mathcal{E}_h$, we denote $\langle u,v\rangle_{e}=\int_{e} uv \mathrm{d} s $ for $u, v \in L^2(e)$ when $d=2$, and $\langle u,v\rangle_{e} = u v $ for two real-value scalars $u$ and $v$ when $d=1$. The discrete volume inner products of $u,v \in L^2(\Omega)$ and $\boldsymbol{w}, \boldsymbol{v} \in \left(L^2(\Omega)\right)^d$ are defined as

\begin{equation}
        (u, v)_{\mathcal{T}_h}=\sum_{K \in \mathcal{T}_h}(u, v)_K, \quad (\boldsymbol{w}, \boldsymbol{v})_{\mathcal{T}_h}=\sum_{K \in \mathcal{T}_h}(\boldsymbol{w}, \boldsymbol{v})_K,
\end{equation}

\noindent
and edge inner product of $u,v \in L^2(\mathcal{E}_h)$ is given by

\begin{equation}
    \langle u, v\rangle_{\partial \mathcal{T}_h}=\sum_{K \in \mathcal{T}_h}\langle u, v\rangle_{\partial K}, \quad \langle u, v\rangle_{D}=\sum_{e \in D}\langle u, v\rangle_{e},
\end{equation}
for any $D \subseteq \mathcal{E}_h$.

The HDG scheme for equations \eqref{FirstOrderGoverningEquations} seeks approximations $(\mathbf{E}_h, \psi_h, \widehat{\psi}_h)\in \boldsymbol{V}_h^p\times W_h^p\times M_h^p$, $(\mathbf{J}_{n,h}, n_h, \widehat{n}_h)\in \boldsymbol{V}_h^p\times W_h^p\times M_h^p$ and $(\mathbf{J}_{p,h}, p_h, \widehat{p}_h)\in \boldsymbol{V}_h^p\times W_h^p\times M_h^p$ such that for all $(\mathbf{K},\phi) \in \boldsymbol{V}_h^p\times W_h^p$,
\begin{equation}
    \label{HDGWeakForm}
    \left\{
    \begin{aligned}
    & \left(\mathbf{K}, \mathbf{E}_h\right)_{\mathcal{T}_h}-\left(\nabla \cdot \mathbf{K}, \psi_h\right)_{\mathcal{T}_h}+\left\langle\mathbf{K} \cdot \mathbf{n}, \widehat{\psi_h}\right\rangle_{\partial \mathcal{T}_h} = 0, \\
    & \left(\mathbf{K}, \mathbf{J}_{n,h}\right)_{\mathcal{T}_h}-\mu_n\left(\mathbf{K}, n_h \mathbf{E}_h\right)_{\mathcal{T}_h}+D_n\left(\nabla \cdot \mathbf{K}, n_h\right)_{\mathcal{T}_h}-D_n\left\langle\mathbf{K} \cdot \mathbf{n}, \widehat{n_h}\right\rangle_{\partial \mathcal{T}_h} = 0, \\
    & \left(\mathbf{K}, \mathbf{J}_{p,h}\right)_{\mathcal{T}_h}-\mu_p\left(\mathbf{K}, p_h \mathbf{E}_h\right)_{\mathcal{T}_h}-D_p\left(\nabla \cdot \mathbf{K}, p_h\right)_{\mathcal{T}_h}+D_p\left\langle\mathbf{K} \cdot \mathbf{n}, \widehat{p_h}\right\rangle_{\partial \mathcal{T}_h} = 0, \\
    & -\lambda^2\left(\boldsymbol{\nabla} \phi, \mathbf{E}_h\right)_{\mathcal{T}_h}+\left(\phi, n_h\right)_{\mathcal{T}_h}-\left(\phi, p_h\right)_{\mathcal{T}_h}-\left(\phi, N\right)_{\mathcal{T}_h}+\lambda^2\left\langle\phi, \widehat{\mathbf{E}_h}\cdot \mathbf{n}\right\rangle_{\partial \mathcal{T}_h} = 0, \\
    & \left(\boldsymbol{\nabla} \phi, \mathbf{J}_{n,h}\right)_{\mathcal{T}_h}-\left\langle \phi, \widehat{\mathbf{J}_{n,h}} \cdot \mathbf{n}\right\rangle_{\partial \mathcal{T}_h}+\left(\phi, R_n\right)_{\mathcal{T}_h} = 0, \\
    & -\left(\boldsymbol{\nabla} \phi, \mathbf{J}_{p,h}\right)_{\mathcal{T}_h}+\left\langle \phi, \widehat{\mathbf{J}_{p,h}} \cdot \mathbf{n}\right\rangle_{\partial \mathcal{T}_h}+\left(\phi, R_p\right)_{\mathcal{T}_h} = 0.
    \end{aligned}
    \right.
\end{equation}
\noindent
Numerical fluxes $\widehat{\mathbf{E}_h}, \widehat{\mathbf{J}_{n,h}}, \widehat{\mathbf{J}_{p,h}}$ are given by
\begin{equation}
    \label{HDGFluxes}
    \left\{
    \begin{aligned}
        \widehat{\mathbf{E}_h} & =\mathbf{E}_h+\tau_\psi\left(\psi_h-\widehat{\psi_h}\right) \mathbf{n}, \\
        \widehat{\mathbf{J}_{n,h}} & =\mathbf{J}_{n,h}+\tau_n\left(n_h-\widehat{n_h}\right) \mathbf{n}, \\
        \widehat{\mathbf{J}_{p,h}} & =\mathbf{J}_{p,h}+\tau_p\left(p_h-\widehat{p_h}\right) \mathbf{n},
    \end{aligned}
    \right.
\end{equation}
where $\tau_\psi,\tau_n,\tau_p$ are the local stabilization parameters \citep{CockburnLinear2009}. For the Poisson equation, we simply define $\tau_\psi=1$. For the convection-diffusion equations, we separate convection and diffusion terms in a similar way as \citep{CockburnLinear2009}
\begin{equation}
    \label{StabilizationParameters}
    \begin{aligned}
    &\tau_n=\tau_{n,{\rm conv}}+\tau_{n,{\rm diff}}, &\tau_{n,{\rm conv}}=\mu_n\left|\widehat{\mathbf{E}_h}\right|_{\partial \mathcal{T}_h}, & &\tau_{n,{\rm diff}}=\frac{D_n}{\ell_n},\\
    &\tau_p=\tau_{p,{\rm conv}}+\tau_{p,{\rm diff}}, &\tau_{p,{\rm conv}}=\mu_p\left|\widehat{\mathbf{E}_h}\right|_{\partial \mathcal{T}_h}, & &\tau_{p,{\rm diff}}=\frac{D_p}{\ell_p},
    \end{aligned}
\end{equation}
with $\ell_n$ and $\ell_p$ being the diameter of a specific cell. The continuities of the numerical ﬂuxes along the normal directions for the interior edges $\mathcal{E}_h^o$ and the Neumann boundary conditions along $\partial\Omega_{\rm Neumann}$  are enforced by the following transmission conditions
\rtwo{
\begin{equation}
    \label{HDGTransmission}
    \left\{
    \begin{aligned}
    & \left\langle\lambda, \llbracket \widehat{\mathbf{E}_h} \cdot \mathbf{n}  \rrbracket \right\rangle_{\partial \mathcal{E}_h^o\cup\partial\Omega_{\rm Neumann}}=0, \\
    & \left\langle\lambda, \llbracket \widehat{\mathbf{J}_{n,h}} \cdot \mathbf{n} \rrbracket \right\rangle_{\partial \mathcal{E}_h^o\cup\partial\Omega_{\rm Neumann}}=0, \\
    & \left\langle\lambda, \llbracket \widehat{\mathbf{J}_{p,h}} \cdot \mathbf{n} \rrbracket  \right\rangle_{\partial \mathcal{E}_h^o\cup\partial\Omega_{\rm Neumann}}=0,
    \end{aligned}
    \right.
\end{equation}
\noindent
where $\lambda\in M_h^p$, and $\llbracket \cdot \rrbracket $ denotes the jump to be defined as follows. For an interior edge $e = \partial K^+ \cap \partial K^- \in \mathcal{E}_h^o$, we denote $\mathbf{q}^{\pm}$ as the trace of a given vector $\mathbf{q}$ from $K^{\pm}$ to $e$, and $\mathbf{n}^{\pm}$ as the unit outward normal vector on $\partial K^{\pm}$. Consequently, the jump in the normal component of $\mathbf{q}$ across $e$ is given by $\llbracket \mathbf{q} \cdot \mathbf{n} \rrbracket = \mathbf{q}^+ \cdot \mathbf{n}^+ + \mathbf{q}^- \cdot \mathbf{n}^-$. For a boundary edge $e \in \mathcal{E}^\partial_h$, $\mathbf{q} \cdot \mathbf{n}$ is single-value and therefore $\llbracket \mathbf{q} \cdot \mathbf{n} \rrbracket = \mathbf{q} \cdot \mathbf{n}$. }

The Dirichlet boundary conditions are imposed to trace variables by
\begin{equation}
    \label{HDG_Dirichlet_BC}
    \left\{
    \begin{aligned}
    & \left\langle\lambda, \widehat{\psi_h}-\psi|_{\partial \Omega_{\rm Dirichlet}}\right\rangle_{\partial \Omega_{\rm Dirichlet}}=0, \\
    & \left\langle\lambda, \widehat{n_h}-n|_{\partial \Omega_{\rm Dirichlet}}\right\rangle_{\partial \Omega_{\rm Dirichlet}}=0, \\
    & \left\langle\lambda, \widehat{p_h}-p|_{\partial \Omega_{\rm Dirichlet}}\right\rangle_{\partial \Omega_{\rm Dirichlet}}=0.
    \end{aligned}
    \right.
\end{equation}
Now, \eqref{HDGWeakForm}, \eqref{HDGFluxes}, \eqref{HDGTransmission} and \eqref{HDG_Dirichlet_BC} together with boundary conditions \eqref{NeumannBoundary} and \eqref{DirichletBoundaryRes} complete the {\bf HDG scheme} for the time-independent PNP equations \eqref{FirstOrderGoverningEquations}.

\subsection{The $hp$-adaptivity of the HDG scheme}

\begin{figure}[!ht]
\centering
\centering
\includegraphics[width=0.36\textwidth]{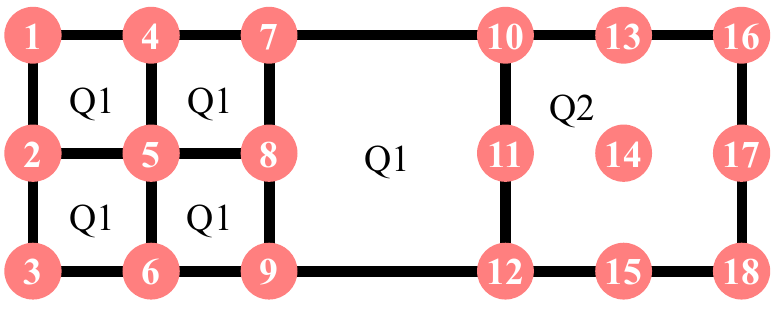}
\caption{Illustration of $hp$-adaptivity on a quadrilateral mesh.}
\label{fig:hp_adaptive}
\end{figure}

One of the biggest advantages of discontinuous Galerkin schemes is that they are naturally suited for $hp$-adaptivity, where elements of different refinement levels and polynomial degrees are solved simultaneously to get accurate results with less unknowns. Fig. \ref{fig:hp_adaptive} shows how HDG scheme handles $hp$-adaptivity on a 2D quadrilateral mesh, where red circles denote degree of freedoms (DoFs) of the skeleton variables $\widehat{\psi_h},\widehat{p_h},\widehat{n_h}$. DoF 8 and DoF 11 are constrained DoFs, originating from $h$-adaptivity and $p$-adaptivity, respectively. By enforcing the constraints $u_8=\frac{u_7+u_9}{2},u_{11}=\frac{u_{10}+u_{12}}{2}$, all DoFs can be uniquely solved.

\subsection{The Newton method and static condensation of HDG scheme}

After discretization, \eqref{HDGWeakForm} and \eqref{HDGTransmission} can  be solved with the Newton method. In each Newton iteration, a global linear system is solved for solution update

\begin{equation}
    \label{HDGinMatrixForm}
    \begin{bmatrix}
        \mathbf{A} & \mathbf{B} \\
        \mathbf{C} & \mathbf{D} \\
    \end{bmatrix}
    \begin{bmatrix}
        \Delta\mathbf{U} \\ \Delta\boldsymbol{\Lambda}
    \end{bmatrix}
    =
    \begin{bmatrix}
        \mathbf{F} \\ \mathbf{G} 
    \end{bmatrix},
\end{equation}

\noindent
where $\Delta\mathbf{U}$ represents the Newton step of solutions $\mathbf{E}_h,\mathbf{J}_{n,h},\mathbf{J}_{p,h},\psi_h,n_h,p_h$, and $\Delta\boldsymbol{\Lambda}$ represents the Newton step of $\widehat{\psi_h},\widehat{n_h},\widehat{p_h}$.

Equations in \eqref{HDGWeakForm} can be solved on each cell independently as $\mathbf{E}_h,\mathbf{J}_{n,h},\mathbf{J}_{p,h}$ live on discontinuous approximation space $\boldsymbol{V}_h^p$ and $\psi_h,n_h,p_h$ live on discontinuous space $W_h^p$. In other words, the matrix $\mathbf{A}$ is block diagonal after a permutation of $\Delta\mathbf{U}$, with each block size at the order of $p^d$. We can therefore apply the static condensation technique \citep{Cockburn2016Static} and reduce equation \eqref{HDGinMatrixForm} to a global linear system with only $\Delta\boldsymbol{\Lambda}$ as variable

\begin{equation}
    \label{HDGStaticCondensation}
    (\mathbf{D}-\mathbf{C}\mathbf{A}^{-1}\mathbf{B})\Delta\boldsymbol{\Lambda}=\mathbf{G}-\mathbf{C}\mathbf{A}^{-1}\mathbf{F}.
\end{equation}
After solving the global linear system \eqref{HDGStaticCondensation} for $\Delta\boldsymbol{\Lambda}$, solution $\Delta\mathbf{U}$ can be obtained locally on each cell by 

\begin{equation}
    \label{HDGLocalReconstruct}
    \mathbf{A}\Delta\mathbf{U}=\mathbf{F}-\mathbf{B}\Delta\boldsymbol{\Lambda},
\end{equation}

\noindent
which is ideally suited for parallel computation.

For our numerical simulation, two convergence criteria are utilized for the Newton method:
\begin{enumerate}
    \item Absolute tolerance: the maximum residual of all the equations \eqref{HDGWeakForm} and \eqref{HDGFluxes} are both below a given tolerance $\varepsilon_{tol}$, e.g. $\varepsilon_{tol}=10^{-8}$,
    
    \item Relative tolerance: the maximum relative Newton step sizes $\max\left(\frac{|\delta u|}{|u| + \varepsilon_{abs}}\right)$ are below a given threshold $\varepsilon_{rel}$, e.g. $\varepsilon_{rel}=10^{-6}$, where $u$ represents scalar variables $\psi_h, n_h, p_h, \widehat{\psi_h}, \widehat{n_h}, \widehat{p_h}$ and $\varepsilon_{abs}=1$ for $\psi_h, \widehat{\psi_h}$ and $\varepsilon_{abs}=n_{\rm ie}$ for $n
    _h,p_h,\widehat{n_h}, \widehat{p_h}$.
\end{enumerate}
The Newton iteration is considered to converge  if one of the above two criteria is met. This completes the Newton method for solving coupled nonlinear semiconductor device equations  by the HDG scheme.

\subsection{Local post-processing of the HDG scheme}

The numerical solution $\psi_h, n_h,p_h$ of the HDG scheme can achieve super-convergence by cell-wise post-processing \cite{CockburnNonlinear2009}. Unlike the approach in \cite{CockburnNonlinear2009} where all variables are post-processed, we focus solely on post-processing the scalar variables for simplicity, as post-processing is not the primary focus of this paper. The post-processed $\psi_h^*, n_h^*,p_h^*$ can reach $k+2$ order accuracy in the $L^2$-norm when $\psi_h, n_h, p_h$ are approximated by polynomials of order $k$. These scalar variables $\psi_h^*, n_h^*,p_h^* \in W_h^{p+1}$ are obtained by solving the local problem on cell $K$
\begin{equation}
    \label{PostProcessing}
    \left\{
        \begin{aligned}
            &\left(1, \psi_h^*\right)_K-\left(1, \psi_h\right)_K=0,\\
            &\left(1, n_h^*\right)_K-\left(1, n_h\right)_K=0,\\
            &\left(1, p_h^*\right)_K-\left(1, p_h\right)_K=0,\\
            &\left(\nabla w_h^*, \psi^*_h\right)_K + \left(\nabla w_h^*, \mathbf{E}_h\right)_K=0,\\
            &\left(\nabla w_h^*, \mathbf{J}_{n,h}\right)_K-\mu_n\left(\nabla w_h^*, n^*_h \mathbf{E}_h\right)_K-D_n\left(\nabla w_h^*,\nabla n^*_h\right)_K=0,\\
            &\left(\nabla w_h^*, \mathbf{J}_{p,h}\right)_K-\mu_p\left(\nabla w_h^*, p^*_h \mathbf{E}_h\right)_K+D_p\left(\nabla w_h^*,\nabla p^*_h\right)_K=0,
        \end{aligned}
    \right.
    \qquad \forall w^*_h \in W_h^{p+1}(K).
\end{equation}

\section{HDG scheme with harmonic averaging technique}
\label{sec:HA}
\begin{figure}[!ht]
\centering
\includegraphics[width=0.5\textwidth]{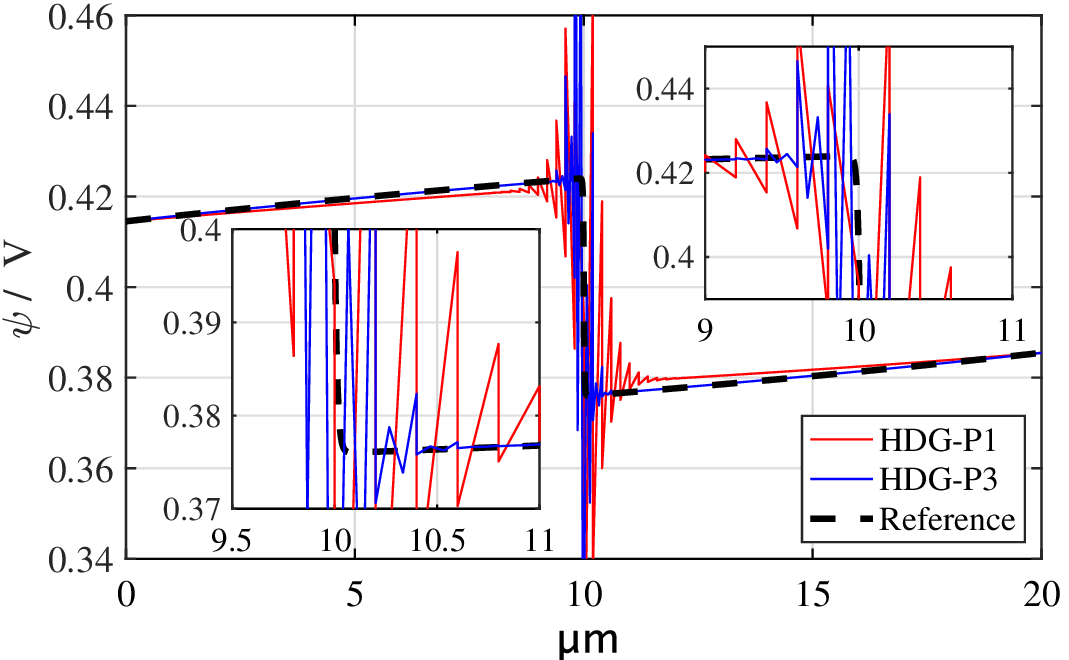}
\caption{The oscillating numerical solution of $\psi$ using the conventional HDG scheme with 100 cells. The settings are taken from test case 3 in subsection \ref{sec:convergeHA}. Reference solution is obtained with FVSG scheme using 10000 cells.}
\label{fig:HDG_oscillation_illustration}
\end{figure}

In semiconductor devices, physical quantities $\psi, n$ and $p$ exhibit sharp jumps near abrupt PN-junctions and form thin interior layers that pose significant challenges to numerical stability. The steep gradient of $\psi$ leads to a large convection speeds $\mu_n \mathbf{E}$ and $\mu_p\mathbf{E}$ in carrier convection-diffusion equations, resulting in non-physical oscillations for conventional HDG scheme \eqref{HDGWeakForm}, \eqref{HDGFluxes} and \eqref{HDGTransmission}. To illustrate these oscillations, we  present the results of the conventional HDG scheme applied to a one-dimensional semiconductor with a single PN-junction in Fig. \ref{fig:HDG_oscillation_illustration}, where non-physical oscillations are clearly observed. Therefore, special treatments have to be imposed to enhance the convergence and stability of the HDG scheme. Inspired by Zhang \textit{et al.} \cite{Zhang2022IAFEM}, we introduce the harmonic averaging technique to improve numerical stability.

\subsection{Harmonic averaging technique and the HDG discretization}

In the HDG discretization \eqref{HDGWeakForm}, convection terms and diffusion terms are treated independently in the auxiliary equations, 
\begin{equation}
    \label{HDGAuxiliaryConvectionDiffusion}
    \begin{aligned}
    & 0=\left(\mathbf{K}, \mathbf{J}_{n,h}\right)_{\mathcal{T}_h}-\mu_n\left(\mathbf{K}, n_h \mathbf{E}_h\right)_{\mathcal{T}_h}+D_n\left(\nabla \cdot \mathbf{K}, n_h\right)_{\mathcal{T}_h}-D_n\left\langle\mathbf{K} \cdot \mathbf{n}, \widehat{n_h}\right\rangle_{\partial \mathcal{T}_h}, \\
& 0=\left(\mathbf{K}, \mathbf{J}_{p,h}\right)_{\mathcal{T}_h}-\mu_p\left(\mathbf{K}, p_h \mathbf{E}_h\right)_{\mathcal{T}_h}-D_p\left(\nabla \cdot \mathbf{K}, p_h\right)_{\mathcal{T}_h}+D_p\left\langle\mathbf{K} \cdot \mathbf{n}, \widehat{p_h}\right\rangle_{\partial \mathcal{T}_h}.
    \end{aligned}
\end{equation}
To obtain a robust numerical discretization in the convection-dominated regime, we combine the convection and diffusion terms in \eqref{HDGAuxiliaryConvectionDiffusion} together as
\begin{equation}
    \label{HDGAuxiliaryReIntegrate1}
    \begin{aligned}
     0&=\left(\mathbf{K}, \mathbf{J}_{n,h}\right)_{\mathcal{T}_h}
    +D_n\left(\mathbf{K},  \frac{1}{V_n}n_h \nabla \psi_h -\nabla n_h\right)_{\mathcal{T}_h}-D_n\left\langle\mathbf{K} \cdot \mathbf{n}, \widehat{n_h}-n_h\right\rangle_{\partial \mathcal{T}_h} \\
    &=\left(\mathbf{K}, \mathbf{J}_{n,h}\right)_{\mathcal{T}_h}
    +D_n\left(\mathbf{K}, -\exp\left(\frac{\psi_h}{V_n}\right)\nabla\left(n_h\exp\left(-\frac{\psi_h}{V_n}\right)\right)\right)_{\mathcal{T}_h}
    -D_n\left\langle\mathbf{K} \cdot \mathbf{n}, \widehat{n_h}-n_h\right\rangle_{\partial \mathcal{T}_h}, \\
     0&=\left(\mathbf{K}, \mathbf{J}_{p,h}\right)_{\mathcal{T}_h}
    +D_p\left(\mathbf{K}, \frac{1}{V_p}p_h \nabla \psi_h + \nabla p_h \right)_{\mathcal{T}_h}+D_p\left\langle\mathbf{K} \cdot \mathbf{n}, \widehat{p_h} - p_h\right\rangle_{\partial \mathcal{T}_h} \\
    & =\left(\mathbf{K}, \mathbf{J}_{p,h}\right)_{\mathcal{T}_h}
   +D_p\left(\mathbf{K}, \exp\left(-\frac{\psi_h}{V_p}\right)\nabla\left(p_h\exp\left(\frac{\psi_h}{V_p}\right)\right)\right)_{\mathcal{T}_h}
    +D_p\left\langle\mathbf{K} \cdot \mathbf{n}, \widehat{p_h} - p_h\right\rangle_{\partial \mathcal{T}_h},
    \end{aligned}
\end{equation}
where $V_n=\frac{D_n}{\mu_n}$, $V_p=\frac{D_p}{\mu_p}$,
and then we take a similar approach as \citep{Zhang2022IAFEM} to approximate the exponentially fitted concentration $n_h\exp\left(-\frac{\psi_h}{V_n}\right), p_h\exp\left(\frac{\psi_h}{V_p}\right)$. In the remainder of this subsection, we will address the equation of $\mathbf{J}_{n,h}$ in detail to illustrate our idea, while the equation of $\mathbf{J}_{p,h}$ can be handled in a similar manner.

 We notice that \eqref{HDGAuxiliaryReIntegrate1} is equivalent to rewrite the linear drift-diffusion operator $-D_n\nabla^2n_h+\mu_n\nabla(n_h\nabla\psi_h)$ (w.r.t. $n_h$) in a self adjoint form $-D_n\nabla(e^{\psi_h/V_n}\nabla(n_he^{-\psi_h/V_n}))$ (w.r.t the new variable $n_he^{-\psi_h/V_n}$). 

 To mimic the SG discretization, we consider linear elements, i.e., $p=1$ is chosen in \eqref{DGSpaces} and \eqref{TraceSpace} with the corresponding finite element spaces $\boldsymbol{V}_h^1$, ${W}_h^1$ and $M_h^1$. 
For $K\in\mathcal{T}_h$, let $K_i$ ($i=1,2,\ldots,N_K$) be the element vertices, and denote the Lagrange basis functions of the spaces $\boldsymbol{V}_h^1$ and ${W}_h^1$ by $\mathbf{K}_{i,m}$ ($m=1,\ldots,d$) and $s_j$, respectively.
More specifically, at vertex $K_i$, only $m$-th entry of $\mathbf{K}_{i,m}$ is non-zero, and $\mathbf{K}_{i,m}$ is equal to the zero vector at all other vertices. Similarly, $s_j=1$ at vertex $K_j$ and $s_j=0$ at all other vertices.

With the help of the above Lagrange basis functions, any function $\mathbf{K}\in\boldsymbol{V}_h^1$ and $\phi \in {W}_h^1$ can be expanded in terms of their nodal values at the vertices as follows:
\begin{equation}\label{eq:nodalexpansion}
\mathbf{K}=\sum\limits_{K\in\mathcal{T}_h}\sum\limits_{i=1}^{N_K}\sum\limits_{m=1}^dK_{i,m}\mathbf{K}_{i,m}, \quad \mathbf{\phi}=\sum\limits_{K\in\mathcal{T}_h}\sum\limits_{j=1}^{N_K} \phi_{j} s_{j}.
\end{equation}
In any element $K$, we denote the nodal values of $\psi_h$, $p_h$ and $n_h$ at the vertices $K_i$ as $\psi_i^K$, $p_i^K$ and $n_i^K$, respectively.

Plugging the nodal expansion \eqref{eq:nodalexpansion} into the volume integral of the convection-diffusion term in \eqref{HDGAuxiliaryReIntegrate1}, we obtain the following approximation
\begin{equation}\label{HAapp1}
    D_n\left(\mathbf{K}, -\exp\left(\frac{\psi_h}{V_n}\right)\nabla\left(n_h\exp\left(-\frac{\psi_h}{V_n}\right)\right)\right)_{\mathcal{T}_h}
    \approx -D_n\sum_{K \in \mathcal{T}_h} \sum_{i=1}^{N_K} \sum_{j=1}^{N_K}\sum\limits_{m=1}^d  K_{i,m} n_j^K\exp\left(-\frac{\psi_j^K}{V_n}\right) \int_K\exp\left(\frac{\psi_h}{V_n}\right)\mathbf{K}_{i,m}\cdot \nabla s_j d\mathbf{x}, \ \  
\end{equation}
where the integrand $n_h\exp\left(-\frac{\psi_h}{V_n}\right)$ is approximated by its Lagrange interpolation as $n_h\exp\left(-\frac{\psi_h}{V_n}\right)\approx\sum_{j=1}^{N_K}n_j^K\exp\left(-\frac{\psi_j^K}{V_n}\right)s_j$. For shape functions $s_j$, it is easy to verify that $\nabla s_j=0-\sum_{i\neq j,i=1}^{N_K} \nabla s_i$. Therefore, if we denote
\begin{equation}\label{integrala}
 e^{\psi_h,K}_{ij,m}=\int_K \exp\left(\frac{\psi_h}{V_n}\right)\mathbf{K}_{i,m}\cdot \nabla s_j d\mathbf{x},
\end{equation} 

it holds that
\begin{equation}
    \label{eij}
    e^{\psi_h,K}_{ii,m}=0-\sum_{j\neq i,j=1}^{N_k} e^{\psi_h,K}_{ij,m}.
\end{equation}

Making use of \eqref{eij}, we can reformulate \eqref{HAapp1} as
\begin{equation}\label{HAapp2}
    \begin{aligned}
    &D_n\left(\mathbf{K}, -\exp\left(\frac{\psi_h}{V_n}\right)\nabla\left(n_h\exp\left(-\frac{\psi_h}{V_n}\right)\right)\right)_{\mathcal{T}_h}\\
    &\approx   D_n\sum_{K \in \mathcal{T}_h} \sum_{i=1}^{N_K} \sum_{j=1,j\neq i}^{N_K}\sum\limits_{m=1}^d  K_{i,m} \left(n_i^K\exp\left(-\frac{\psi_i^K}{V_n}\right)-n_j^K\exp\left(-\frac{\psi_j^K}{V_n}\right)\right)e^{\psi_h,K}_{ij,m}.
    \end{aligned}
\end{equation}

For $e^{\psi_h,K}_{ij,m}$ in \eqref{integrala}, one can approximate $\exp\left(\frac{\psi_h}{V_n}\right)$ by a constant function $\frac{1}{2}\left(\exp\left(\frac{\psi_i^K}{V_n}\right)+\exp\left(\frac{\psi_j^K}{V_n}\right)\right)$ as
\begin{equation}\label{2ndapp}
 e^{\psi_h,K}_{ij,m}=\int_K \exp\left(\frac{\psi_h}{V_n}\right)\mathbf{K}_{i,m}\cdot \nabla s_j dK\approx\frac{1}{2}\left(\exp\left(\frac{\psi_i^K}{V_n}\right)+\exp\left(\frac{\psi_j^K}{V_n}\right)\right) \int_{K}\mathbf{K}_{i,m}\cdot \nabla s_j \,dK.
\end{equation}
Denoting
\begin{equation}
e_{ij,m}^K= \int_{K}\mathbf{K}_{i,m}\cdot \nabla s_j \,dK,
\end{equation}
one arrives at an approximation of \eqref{HAapp2} as
\begin{equation}
    \label{HDGAuxiliaryReIntegrate2}
    \begin{aligned}
   & D_n\left(\mathbf{K}, -\exp\left(\frac{\psi_h}{V_n}\right)\nabla\left(n_h\exp\left(-\frac{\psi_h}{V_n}\right)\right)\right)_{\mathcal{T}_h}\\
    &\approx D_n\sum_{K \in \mathcal{T}_h} \sum_{i=1}^{N_K} \sum_{j=1,j\neq i}^{N_K}\sum\limits_{m=1}^d \frac{1}{2} \left( \exp\left(\frac{\psi_i^K}{V_n}\right)+\exp\left(\frac{\psi_j^K}{V_n}\right) \right) K_{i,m} \left(n_i^K\exp\left(-\frac{\psi_i^K}{V_n}\right)-n_j^K\exp\left(-\frac{\psi_j^K}{V_n}\right)\right)e^{K}_{ij,m}.
    \end{aligned}
\end{equation}
Unfortunately, the above \eqref{HDGAuxiliaryReIntegrate2} may still suffer from non-physical oscillations near junctions where $\psi_h$ exhibits huge jump. Therefore, we propose the following harmonic average $E_{ij}\left(\frac{\psi_h}{V_n}\right)$ to replace  $\frac{1}{2}\left(\exp\left(\frac{\psi_i^K}{V_n}\right)+\exp\left(\frac{\psi_j^K}{V_n}\right)\right)$,
\begin{equation}
    \label{HDGAuxiliaryHarmonicAverageTerms}
    \begin{aligned}
        E_{ij}\left( \frac{\psi_h}{V_n}\right)=\frac{-\frac{\psi_i^K-\psi_j^K}{V_n}}{\exp\left(-\frac{\psi_i^K}{V_n}\right)-\exp\left(-\frac{\psi_j^K}{V_n}\right)} =\left(\frac{\int_{\epsilon_{ij}}\exp\left(-\frac{\psi_{h,\epsilon_{ij}}}{V_n}\right)ds}{L_{ij}}\right)^{-1},
   \end{aligned}
\end{equation}
\noindent
where  $\epsilon_{ij}$ is the line segment connecting vertices $K_i$ and $K_j$,  $\psi_{h,\epsilon_{ij}}$ is a linear function on $\epsilon_{ij}$ with boundary values $\psi_{h,\epsilon_{ij}}|_{K_i}=\psi_i^K$ and $\psi_{h,\epsilon_{ij}}|_{K_j}=\psi_j^K$, $L_{ij}$ is the length of $\epsilon_{ij}$. It is clear from the line integral form in \eqref{HDGAuxiliaryHarmonicAverageTerms} that $E_{ij}\left( \frac{\psi_h}{V_n}\right)$ is a harmonic average.

Substituting \eqref{HDGAuxiliaryHarmonicAverageTerms} into \eqref{HDGAuxiliaryReIntegrate2}  to replace $\frac{1}{2}\left(\exp\left(\frac{\psi_i^K}{V_n}\right)+\exp\left(\frac{\psi_j^K}{V_n}\right)\right)$, we get our final approximation for the convection-diffusion term of $\mathbf{J}_{n,h}$ in \eqref{HDGAuxiliaryReIntegrate1}: 
\begin{align}
    D_n\left(\mathbf{K}, -\exp\left(\frac{\psi_h}{V_n}\right)\nabla\left(n_h\exp\left(-\frac{\psi_h}{V_n}\right)\right)\right)_{\mathcal{T}_h}
    &\approx D_n\sum_{K \in \mathcal{T}_h} \sum_{i=1}^{N_K} \sum_{j=1,j\neq i}^{N_K}\sum_{m=1}^d \left(B\left(\frac{\psi_i^K-\psi_j^K}{V_n}\right)n_i^K - B\left(\frac{\psi_j^K-\psi_i^K}{V_n}\right)n_j^K  \right) e^K_{ij,m}K_{i,m} \notag \\
    & := a_h(\mathbf{K}, n_h,\psi_h), \label{HDGAuxiliaryReIntegrateBernoulli}
\end{align}
    where $B(x)=\frac{x}{\exp(x)-1}$ is the Bernoulli function,
 and we shall replace  $D_n\left(\mathbf{K}, -\exp\left(\frac{\psi_h}{V_n}\right)\nabla\left(n_h\exp\left(-\frac{\psi_h}{V_n}\right)\right)\right)_{\mathcal{T}_h}$ by its discrete weak form $a_h(\mathbf{K}, n_h,\psi_h)$. The discretization in \eqref{HDGAuxiliaryReIntegrateBernoulli} is similar to the SG approach used for discretizing current densities, and is expected to obtain similarly stable and convergent results for the DD equations \eqref{FirstOrderGoverningEquations}.
   
For the $\mathbf{J}_{p,h}$ part, we adopt the same strategy to obtain
\begin{align}
    D_p\left(\mathbf{K}, \exp\left(-\frac{\psi_h}{V_p}\right)\nabla\left(p_h\exp\left(\frac{\psi_h}{V_p}\right)\right)\right)_{\mathcal{T}_h}
    &\approx D_p\sum_{K \in \mathcal{T}_h} \sum_{i=1}^{N_K} \sum_{j=1,j\neq i}^{N_K}\sum_{m=1}^d \left(B\left(\frac{\psi_i^K-\psi_j^K}{V_p}\right)p_j^K - B\left(\frac{\psi_j^K-\psi_i^K}{V_p}\right)p_i^K  \right) e^K_{ij,m}K_{i,m} \notag\\
    &:=b_h(\mathbf{K},p_h,\psi_h). \label{HDGAuxiliaryReIntegrateBernoullip}
\end{align}

The resulting HDG-HA discretization for solving equations \eqref{FirstOrderGoverningEquations} seeks approximations $(\mathbf{E}_h, \psi_h, \widehat{\psi}_h)\in \boldsymbol{V}_h^1\times W_h^1\times M_h^1$, $(\mathbf{J}_{n,h}, n_h, \widehat{n}_h)\in \boldsymbol{V}_h^1\times W_h^1\times M_h^1$ and $(\mathbf{J}_{p,h}, p_h, \widehat{p}_h)\in \boldsymbol{V}_h^1\times W_h^1\times M_h^1$ such that for all $(\mathbf{K},\phi) \in \boldsymbol{V}_h^1\times W_h^1$,
\begin{equation}
    \label{HDGHA_WeakForm}
    \left\{
    \begin{aligned}
    & \left(\mathbf{K}, \mathbf{E}_h\right)_{\mathcal{T}_h}-\left(\nabla \cdot \mathbf{K}, \psi_h\right)_{\mathcal{T}_h}+\left\langle\mathbf{K}\cdot \mathbf{n}, \widehat{\psi_h}\right\rangle_{\partial \mathcal{T}_h} = 0, \\
    & \left(\mathbf{K}, \mathbf{J}_{n,h}\right)_{\mathcal{T}_h}+a_h(\mathbf{K},n_h,\psi_h)-D_n\left\langle\mathbf{K} \cdot \mathbf{n}, \widehat{n_h}-n_h\right\rangle_{\partial \mathcal{T}_h} = 0, \\
    & \left(\mathbf{K}, \mathbf{J}_{p,h}\right)_{\mathcal{T}_h}+b_h(\mathbf{K},p_h,\psi_h)+D_p\left\langle\mathbf{K} \cdot \mathbf{n}, \widehat{p_h}-p_h\right\rangle_{\partial \mathcal{T}_h} = 0, \\
    & -\lambda^2\left(\boldsymbol{\nabla} \phi, \mathbf{E}_h\right)_{\mathcal{T}_h}+\left(\phi, n_h\right)_{\mathcal{T}_h}-\left(\phi, p_h\right)_{\mathcal{T}_h}-\left(\phi, N\right)_{\mathcal{T}_h}+\lambda^2\left\langle\phi, \widehat{\mathbf{E}_h}\cdot \mathbf{n}\right\rangle_{\partial \mathcal{T}_h} = 0, \\
    & \left(\boldsymbol{\nabla} \phi, \mathbf{J}_{n,h}\right)_{\mathcal{T}_h}-\left\langle \phi, \widehat{\mathbf{J}_{n,h}} \cdot \mathbf{n}\right\rangle_{\partial \mathcal{T}_h}+\left(\phi, R_n\right)_{\mathcal{T}_h} = 0, \\
    & -\left(\boldsymbol{\nabla} \phi, \mathbf{J}_{p,h}\right)_{\mathcal{T}_h}+\left\langle \phi, \widehat{\mathbf{J}_{p,h}} \cdot \mathbf{n}\right\rangle_{\partial \mathcal{T}_h}+\left(\phi, R_p\right)_{\mathcal{T}_h} = 0,
    \end{aligned}
    \right.
\end{equation}
where $a_h$ and $b_h$ are defined in \eqref{HDGAuxiliaryReIntegrateBernoulli} and \eqref{HDGAuxiliaryReIntegrateBernoullip}, respectively. The numerical fluxes, continuity and the boundary conditions are the same as in \eqref{HDGFluxes}, \eqref{HDGTransmission} and \eqref{HDG_Dirichlet_BC}. Since the harmonic averaging technique doesn't modify the transmission condition \eqref{HDGTransmission}, all routines for static condensation and the Newton method used to solve the coupled equations remain unchanged. \rtwo{Furthermore, since the last three equations in \eqref{HDGHA_WeakForm} and the transmission condition \eqref{HDGTransmission} are identical to  those in the conventional HDG scheme, the local conservativity of the conventional HDG scheme, as demonstrated in \cite{CockburnNonlinear2009,CockburnLinear2009}, is preserved in our HDG-HA scheme.}

\subsection{HA-adaptivity and its indicator}\label{sec:ha-adaptivity}

DG schemes, which utilize discontinuous polynomial spaces across cells, are well-suited for $hp$-adaptivity that allows for different numerical discretizations in different cells. Our HA-adaptivity is similar to the $p$-adaptivity. The proposed harmonic averaging technique \eqref{HDGAuxiliaryHarmonicAverageTerms} only changes volume integral within cells \eqref{HDGHA_WeakForm}, the transmission condition and static condensation procedure are the same for conventional HDG cells and HDG-HA cells. Therefore we can use different techniques in different cells to synergize the desired robustness and high-order property. \rtwo{It should be noted that the FVSG scheme also have similar robustness to our HDG-HA scheme, suggesting that a hybridization of FVSG and DG schemes could be considered.} \rone{However, the implementation of FVSG scheme requires a dual Voronoi grid \cite{Zhang2022IAFEM}, making hybridization more complex compared to our HA-adaptivity within the HDG framework under a standard mesh, particularly in higher dimensions.This flexibility in applying different techniques to different cells makes our HDG-HA scheme easier to extend to an $hp$-adaptive version compared to the FVSG scheme.}

Within cells where the doping profile has huge jump, which may induce stability issues, we should use the robust harmonic averaging technique. Whereas in the other cells where doping concentration and $\psi,n,p$ change smoothly, conventional HDG schemes can be utilized to increase numerical accuracy while reducing total degrees of freedoms.

Inspired by the intuition of locating the cells within the junction regions, i.e., where the doping profile and solutions exhibit huge jumps or discontinuities, several physics-based indicators are proposed.
\begin{enumerate}
    \item Norm of the gradient of electrostatic potential $\|\nabla \psi \|_{L^2}$.
    \item Norm of the gradient of net doping $\|\nabla N \|_{L^2}$.
    \item Norm of the gradient of carrier concentration $\|\nabla n \|_{L^2}$.
\end{enumerate}
Their performances will be validated in numerical experiments.

\section{Numerical experiments}
\label{sec:Experiments}

\subsection{Stability of HDG schemes on abrupt semiconductor junctions}\label{sec:convergeHA}

We first compare the numerical stability of conventional HDG schemes of various orders with that of the proposed HDG-HA scheme. For reference, the FVSG scheme \cite{Bank1983FVSG} is also included in the comparison. The 
stability
of these schemes are evaluated by the minimum number of cells required for a Newton solver to converge. Additionally, the converged results are plotted to assess the numerical stability with respect to oscillations. To ensure the positivity of the numerical densities $n_h$ and $p_h$, a simple cut-off post-processing could be applied \cite{TongCai2023}. 

The steady state of a 1D silicon diode within the domain $(0,\SI{20}{\micro \metre})$ is simulated. Dirichlet boundary voltages of $0 \rm V$ and $0.8 \rm V$ are applied at boundaries $x = \SI{0}{\micro\metre}$ and $x =  \SI{20}{\micro \metre}$, respectively. The doping profile is discontinuous with an abrupt junction at the middle
\begin{equation}
    \label{TestCase1Doping}
    N(x)=\left\{
    \begin{array}{lcl}
    C_1, & & {x < \SI{10}{\micro \metre},}\\
    \frac{C_1+C_2}{2}, & & {x = \SI{10}{\micro \metre},}\\
    C_2, & & {x > \SI{10}{\micro \metre}.}
    \end{array} \right.
\end{equation}
Five test cases are investigated by by selecting different values in \eqref{TestCase1Doping}:

\begin{enumerate}
    \item $C_1=10^{17} \rm cm^{-3}, C_2=3 \times 10^{17} \rm cm^{-3}$,
    \item $C_1=10^{15} \rm cm^{-3}, C_2=-10^{15} \rm cm^{-3}$,
    \item $C_1=10^{17} \rm cm^{-3}, C_2=-10^{17} \rm cm^{-3}$,
    \item $C_1=10^{19} \rm cm^{-3}, C_2=-10^{19} \rm cm^{-3}$,
    \item $C_1=10^{21} \rm cm^{-3}, C_2=-10^{21} \rm cm^{-3}$,
\end{enumerate}
which range from a pure N-type high-low junction (test case 1, used in \cite{Liu2007LDG2,Shu1995MixedRKDGSemiconductor}) to a heavily doped PN-junction (test case 5). The jump of net doping concentration gradually increases from case 1 to case 5, and we anticipate that the difficulty of convergence will also progressively increase.

The Newton method is used to solve the coupled nonlinear equations \eqref{GoverningEquations} on uniform meshes. A voltage stepping method that steps from a convergent $0 \rm V$ bias solution to a $0.8 \rm V$ bias one is used. Results on a coarse mesh (100 cells), a medium mesh (1000 cells) and a fine mesh (10000 cells) are compared to demonstrate the convergent ability of each method. Table \ref{tab:ConvergenceTable} summarizes the convergence of different methods, where the letters \textbf{S}, \textbf{O} and \textbf{F} stand for:

\begin{table}[!ht]
\caption{
Stability of different schemes on an abrupt junction given in \eqref{TestCase1Doping}.}
\label{tab:ConvergenceTable}
\begin{tabular}{|c|c|c|ccc|ccc|ccc|ccc|}
\hline
          Scheme& FVSG   & HDG-HA & \multicolumn{3}{c|}{HDG-P0}                                        & \multicolumn{3}{c|}{HDG-P1}                                        & \multicolumn{3}{c|}{HDG-P2}                                        & \multicolumn{3}{c|}{HDG-P3}                                        \\ \hline
           Number of cells& $10^2$ & $10^2$ & \multicolumn{1}{c|}{$10^2$} & \multicolumn{1}{c|}{$10^3$} & $10^4$ & \multicolumn{1}{c|}{$10^2$} & \multicolumn{1}{c|}{$10^3$} & $10^4$ & \multicolumn{1}{c|}{$10^2$} & \multicolumn{1}{c|}{$10^3$} & $10^4$ & \multicolumn{1}{c|}{$10^2$} & \multicolumn{1}{c|}{$10^3$} & $10^4$ \\ \hline
Case 1& \green{S}      & \green{S}      & \multicolumn{1}{c|}{O}      & \multicolumn{1}{c|}{O}      & \green{S}      & \multicolumn{1}{c|}{\green{S}}      & \multicolumn{1}{c|}{\green{S}}      & \green{S}      & \multicolumn{1}{c|}{O}      & \multicolumn{1}{c|}{\green{S}}      & \green{S}      & \multicolumn{1}{c|}{O}      & \multicolumn{1}{c|}{\green{S}}      & \green{S}      \\ \hline
Case 2& \green{S}      & \green{S}      & \multicolumn{1}{c|}{\red{F}}      & \multicolumn{1}{c|}{O}      & \green{S}      & \multicolumn{1}{c|}{O}      & \multicolumn{1}{c|}{\green{S}}      & \green{S}      & \multicolumn{1}{c|}{O}      & \multicolumn{1}{c|}{\green{S}}      & \green{S}      & \multicolumn{1}{c|}{O}      & \multicolumn{1}{c|}{\green{S}}      & \green{S}      \\ \hline
Case 3& \green{S}      & \green{S}      & \multicolumn{1}{c|}{\red{F}}      & \multicolumn{1}{c|}{\red{F}}      & \red{F}      & \multicolumn{1}{c|}{O}      & \multicolumn{1}{c|}{O}      & \green{S}      & \multicolumn{1}{c|}{\red{F}}      & \multicolumn{1}{c|}{\green{S}}      & \green{S}      & \multicolumn{1}{c|}{O}      & \multicolumn{1}{c|}{\green{S}}      & \green{S}      \\ \hline
Case 4& \green{S}      & \green{S}      & \multicolumn{1}{c|}{\red{F}}      & \multicolumn{1}{c|}{\red{F}}      & \red{F}      & \multicolumn{1}{c|}{\red{F}}      & \multicolumn{1}{c|}{\red{F}}      & \green{S}      & \multicolumn{1}{c|}{\red{F}}      & \multicolumn{1}{c|}{\red{F}}      & \red{F}      & \multicolumn{1}{c|}{\red{F}}      & \multicolumn{1}{c|}{O}      & \green{S}      \\ \hline
Case 5& \green{S}      & \green{S}      & \multicolumn{1}{c|}{\red{F}}      & \multicolumn{1}{c|}{\red{F}}      & \red{F}      & \multicolumn{1}{c|}{\red{F}}      & \multicolumn{1}{c|}{\red{F}}      & \red{F}      & \multicolumn{1}{c|}{\red{F}}      & \multicolumn{1}{c|}{\red{F}}      & \red{F}      & \multicolumn{1}{c|}{\red{F}}      & \multicolumn{1}{c|}{\red{F}}      & O      \\ \hline
\end{tabular}
\end{table}

\begin{itemize}
    \item \textbf{S}: Success. The scheme converges to 
    $0.8\rm V$ bias successfully without any non-physical oscillations;
    
    \item \textbf{O}: Oscillate. Although the scheme converges to $0.8\rm V$ bias, apparent non-physical oscillations occur;
    
    \item \textbf{F}: Fail. The scheme fails to converge to $0.8\rm V$ bias, even when employing a small voltage stepping $10^{-6} \rm V$.
\end{itemize}

Table \ref{tab:ConvergenceTable} clearly reveals 
that when the jump of doping concentration becomes larger, the performance of all conventional HDG schemes become worse,
while the FVSG and the HDG-HA schemes still converge without any non-physical oscillations. Additionally, we found that in this experiment, both the FVSG and the HDG-HA schemes can converge to a $0.8\rm V$ bias on all the five test cases with only $1$ or $2$ cells, demonstrating their superior numerical stability.

\begin{figure}[!h]
    \centering
    \includegraphics[width=0.82\textwidth]{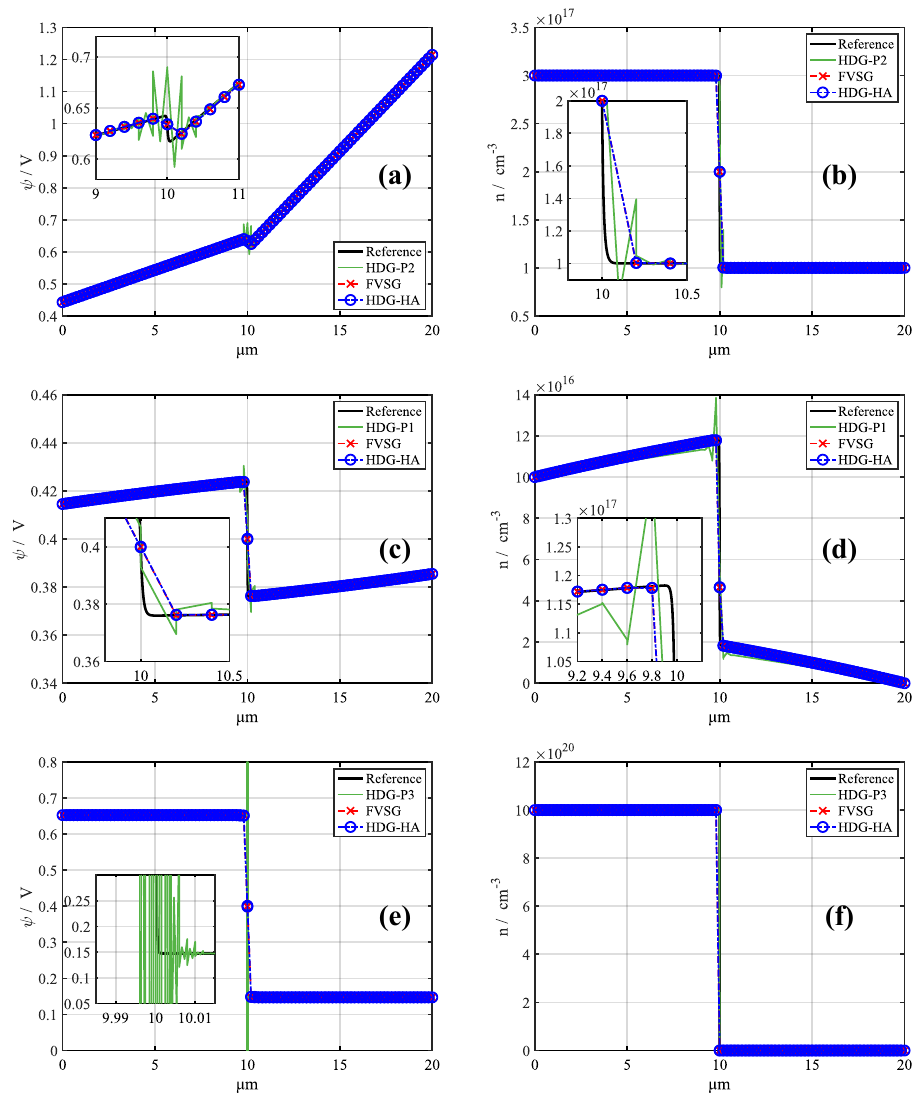}
    \caption{Comparisons of different schemes: solutions of $\psi$ (left column) and solutions of $n$ (right column). Top row: case 1 using $100$ cells. Middle row: case 3 using $100$ cells. Bottom row: case 5 where both HDG-HA and FVSG schemes use $100$ cells, and HDG-P3 scheme used $10000$ cells.
    }
    \label{fig:Exp_1_ResultsCompare}
\end{figure}

Fig. \ref{fig:Exp_1_ResultsCompare} depicts some numerical results in Table \ref{tab:ConvergenceTable}, where settings of the first row ((a) and (b)), second row ((c) and (d)) and third row ((e) and (f)) come from test cases 1, 3 and 5, respectively. The numerical solutions with high-order polynomials (HDG-P2, HDG-P3) 
are linearly interpolated at several points within each cell for the purpose of visualization. As shown in this figure, the solution of the HDG-HA scheme aligns with that of the FVSG scheme, showing no oscillations even at sharp junctions with only 100 cells. In contrast, conventional HDG schemes are prone to numerical oscillations when large gradients are present in the solutions. Subfigures (e) and (f) in Fig. \ref{fig:Exp_1_ResultsCompare} demonstrate that the conventional HDG scheme still oscillates with a fine mesh of $10000$ cells, when applied to a heavily doped PN-junction.

\subsection{Numerical accuracy of various HDG schemes}
As indicated in Table \ref{tab:ConvergenceTable}, 
stability of conventional HDG schemes are challenged by abrupt PN-junctions with huge doping jump.  When the doping profile of devices is smooth enough, conventional HDG schemes show their clear superiority in achieving high-order accuracy.

\begin{figure}[!ht]
\centering
\includegraphics[width=0.32\textwidth]{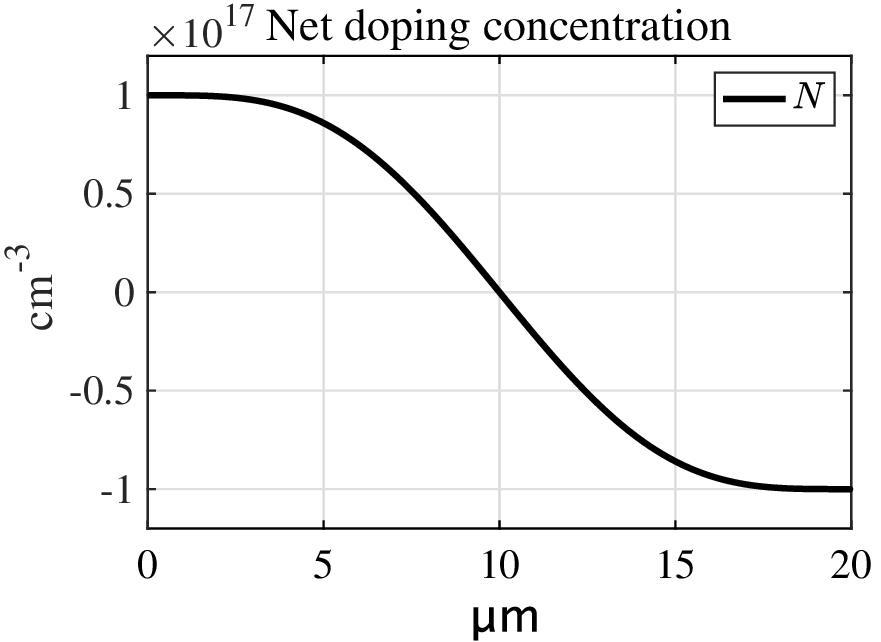}
\includegraphics[width=0.32\textwidth]{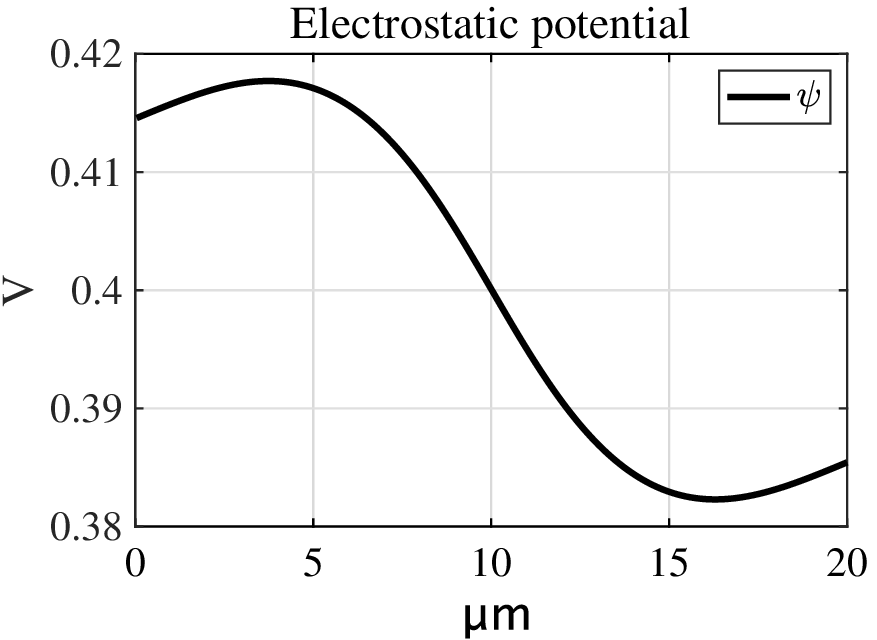}
\includegraphics[width=0.32\textwidth]{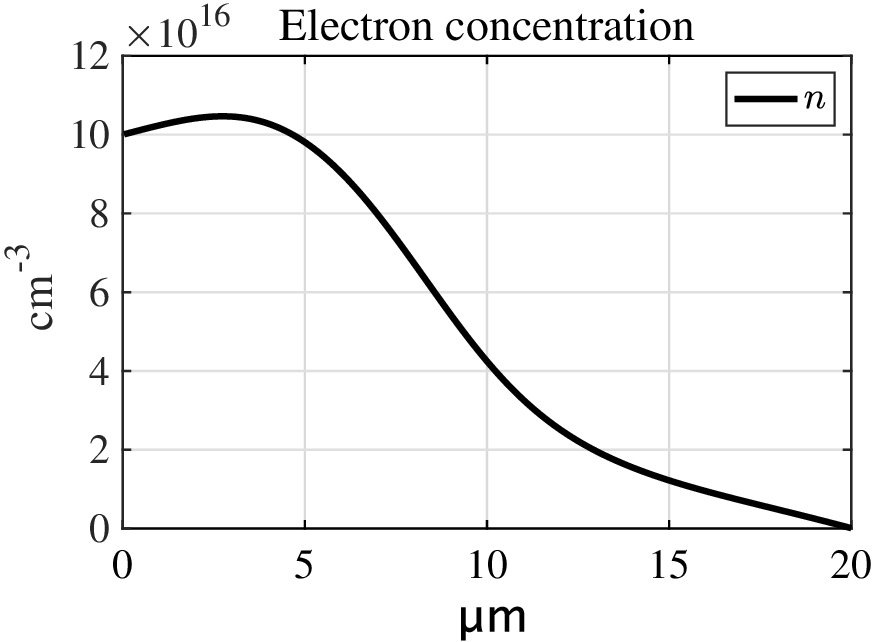}
\caption{The doping and solution profiles of the smooth PN-junction test case, where solutions are computed using the FVSG scheme with $10^6$ cells.
}
\label{fig:Exp_2_SmoothSolution}
\end{figure}

To illustrate this accuracy, we choose a smooth one-dimensional PN diode within the domain $(0,\SI{20}{\micro \metre})$, where the net doping concentration is described by
\begin{equation}
    \label{Exp2SmoothDoping}
   N(x) = \left(1 + 40 \left(\frac{x}{\SI{20}{\micro \metre}}\right)^7-140\left(\frac{x}{\SI{20}{\micro \metre}}\right)^6+168\left(\frac{x}{\SI{20}{\micro \metre}}\right)^5-70\left(\frac{x}{\SI{20}{\micro \metre}}\right)^4\right)\times 10^{17} {\rm cm^{-3}}.
\end{equation}
Voltages of $0 \rm V$ and $0.8 \rm V$ are applied at the boundaries located at $x = \SI{0}{\micro\metre}$ and $x = \SI{20}{\micro\metre}$, respectively. The solution variables $\psi,n,p$ in this example are smooth across the whole computational domain, as depicted in Fig. \ref{fig:Exp_2_SmoothSolution}, and every numerical scheme converges with just a few cells.

\begin{figure}[!ht]
\centering
\includegraphics[width=0.72\textwidth]{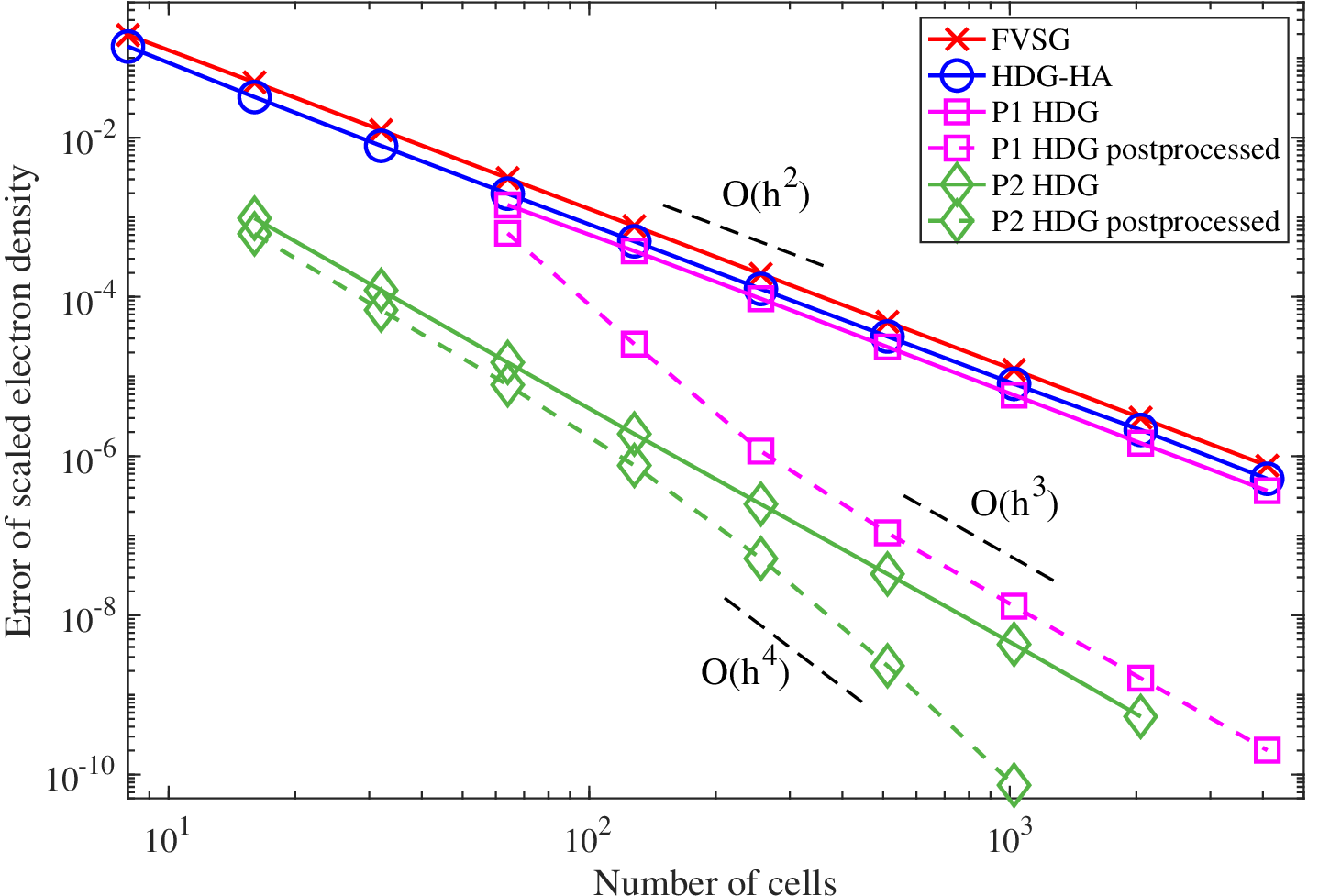}
\caption{$L^\infty$ errors of different schemes
on the smooth PN-junction test case.}
\label{fig:Exp_2_OrderCompare}
\end{figure}

Fig. \ref{fig:Exp_2_OrderCompare} shows the $L^\infty$ errors of the electron density for
the FVSG, HDG-HA and conventional HDG schemes of order 1 and 2. Their $L^2$ errors are similar and are therefore omitted here. From Fig. \ref{fig:Exp_2_OrderCompare}, both the FVSG and the proposed HDG-HA achieve second-order convergence as expected. The errors of the HDG-HA scheme are quite close to those of the FVSG, showing only a minor improvement.
On the contrary, the high-order HDG schemes exhibit significantly more accurate results than both the HDG-HA and FVSG schemes. In addition, the convergence rates can be improved to $O(h^{p+2})$ using polynomials of order $p$, with the help of post-processing. Therefore, it would be advantageous to utilize $hp$-adaptivity to combine the robustness of the HDG-HA scheme and the high-order accuracy of the conventional HDG schemes. In the subsequent discussions, if unspecified, all the HDG solutions (or solutions in the conventional HDG cells) are post-processed to enhance the accuracy.

\subsection{Performance of HA-adaptive indicators}\label{sec:num_ha_indicator}

\rone{
The previous two experiments highlight the robustness of the HDG-HA scheme for solutions with large jumps, as well as the superior accuracy of high-order HDG schemes for smooth solutions. By leveraging the $hp$-adaptivity of HDG schemes, we can effectively combine these two advantages. In this subsection, we will first demonstrate the benefits of this combination using a one-dimensional example and then validate our proposed indicators from subsection \ref{sec:ha-adaptivity} with a two-dimensional example. }     

\begin{figure}[!ht]
\centering
\includegraphics[width=0.36\textwidth]{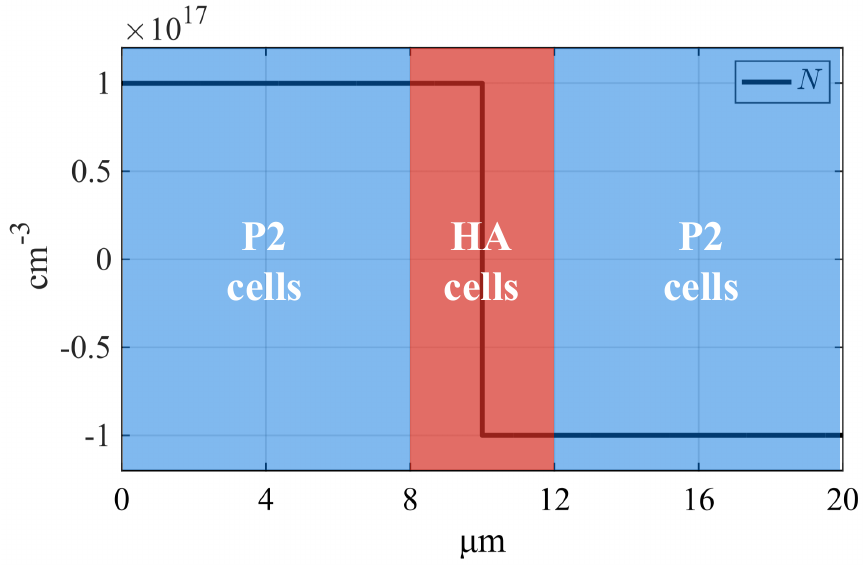}
\caption{\rone{Doping profile and HA+$p$-adaptivity.}}
\label{fig:Supplement_doping}
\end{figure}

\rone{
To validate the effectiveness of our HA-adaptive HDG scheme, we reuse the 1D example, case 3 from subsection \ref{sec:convergeHA} with uniform meshes. The entire domain $(0, \SI{20}{\micro\metre})$ is divided into two distinct regions: the abrupt region $(\SI{8}{\micro\metre}, \SI{12}{\micro\metre})$ and the smooth region $(0, \SI{8}{\micro\metre}) \cup (\SI{12}{\micro\metre}, \SI{20}{\micro\metre})$. The cells in the abrupt region are discretized using our HDG-HA scheme, while the cells in the smooth region adopt the conventional HDG-P2 scheme. This decomposition is illustrated in Fig. \ref{fig:Supplement_doping} using two different colors. We refer to this combined approach as the HDG-HA+P2 scheme in this subsection. The FVSG scheme and the HDG-P2 scheme are also reused for comparison.
}

\begin{figure}[!ht]
\centering
\includegraphics[width=0.6\textwidth]{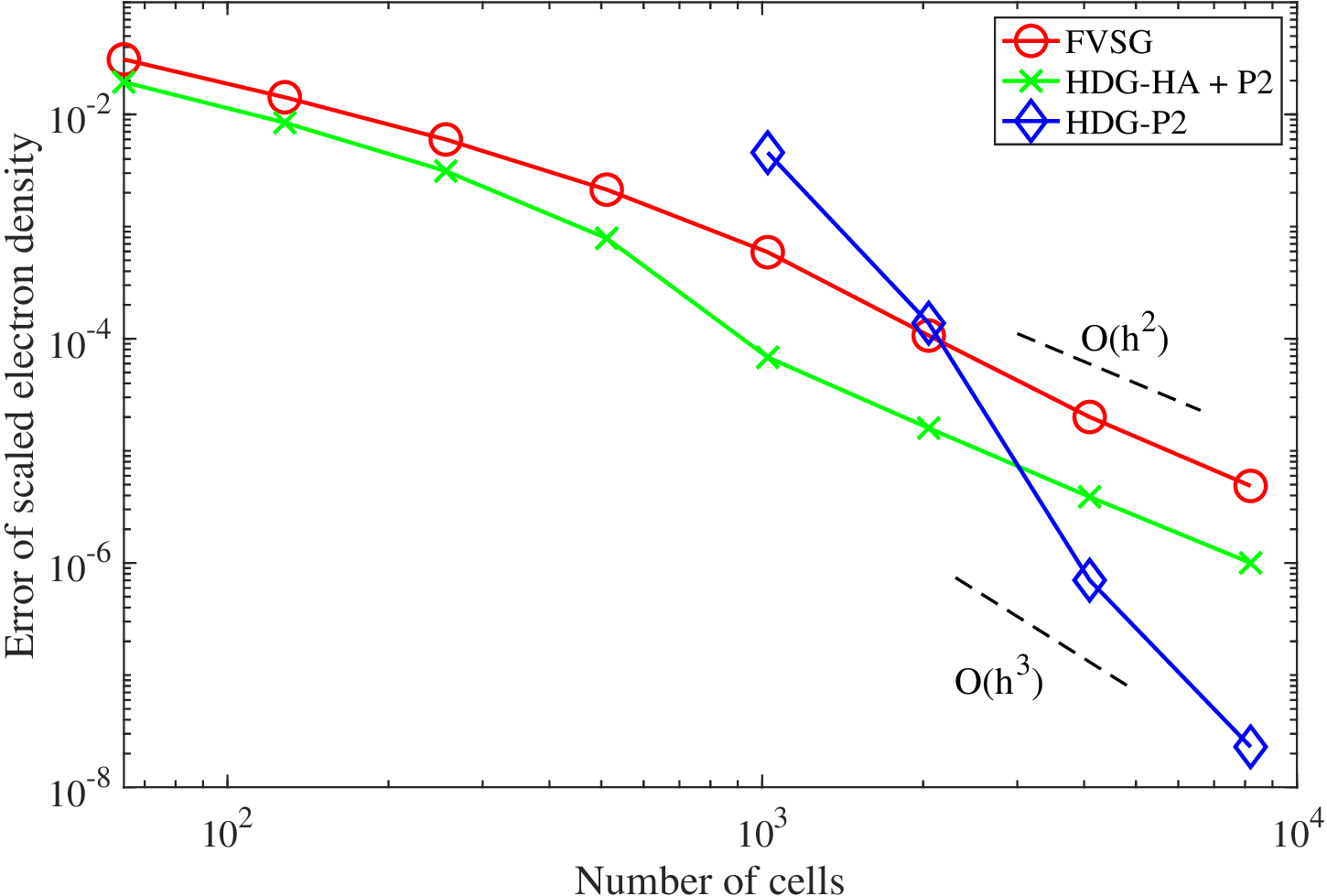}
\caption{\rone{$L^\infty$ errors of different schemes in the smooth region.}}
\label{fig:Supplement_error}
\end{figure}

\rone{
The $L^\infty$ errors of electron density $n$ in the smooth region are shown in Fig. \ref{fig:Supplement_error}. It can be observed that (I) the HDG-HA+P2 scheme reduces errors in smooth region compared to the FVSG scheme, and (II) HDG-HA+P2 scheme gains more stability compared to the HDG-P2 scheme which encounters convergence issues with fewer than $10^3$ cells.
}

\begin{figure}[!ht]
    \centering
    \begin{subfigure}[b]{0.45\textwidth}
        \centering
        \includegraphics[width=0.71\textwidth]{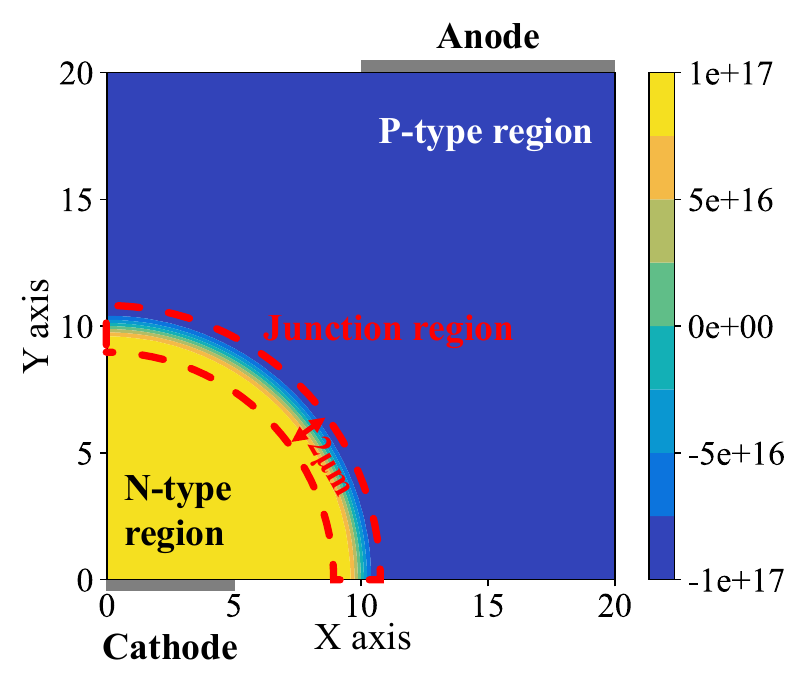}
        \label{fig:Exp_3_2D_case_structure}
    \end{subfigure}
    \begin{subfigure}[b]{0.45\textwidth}
        \centering
        \includegraphics[width=0.8\textwidth]{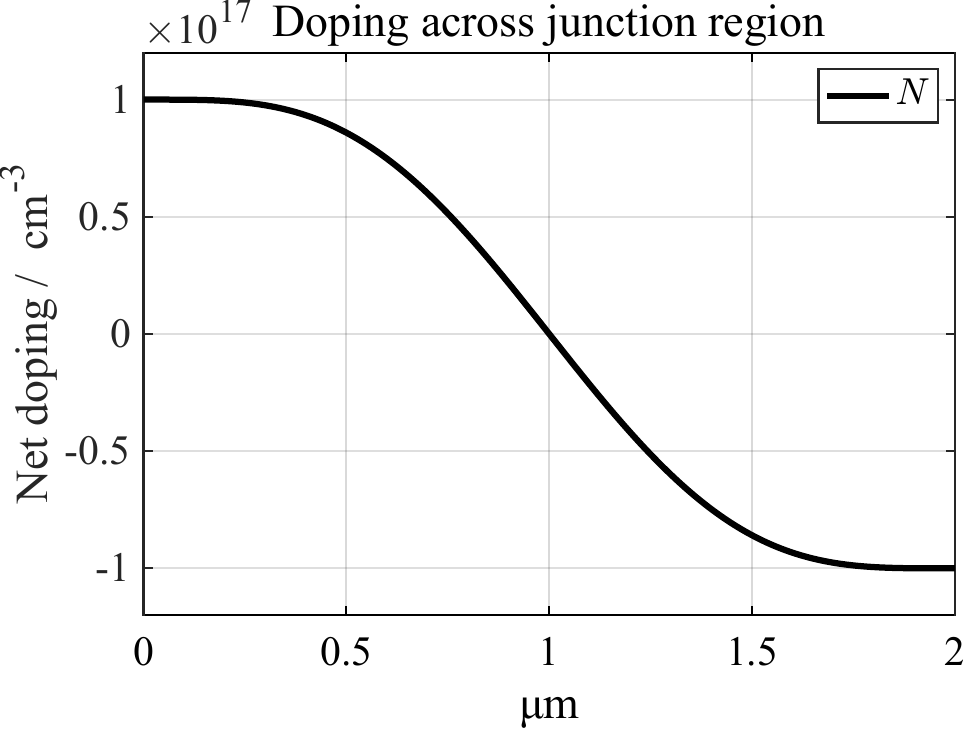}
        \label{fig:Exp_3_2D_case_junction_doping}
    \end{subfigure}
    \caption{A 2D general diode test case with circular junction: structure of the 2D diode (left) and net doping across the junction region (right).}
    \label{fig:Exp_3_2D_case}
\end{figure}

To illustrate the indicators for the HA-adaptive procedure in our HDG-HA scheme, a 2D diode is chosen as the test case, which can serve as a representative for general semiconductor devices. It is square-shaped with a circular PN-junction, as depicted in Fig. \ref{fig:Exp_3_2D_case}. The net doping concentration $N(\textbf{x})$ is given as
\begin{equation}
    \label{Exp3SmoothDoping}
    N(\mathbf{x}) = \hat{N}(r) = \left\{
    \begin{array}{ll}
         10^{17} {\rm cm^{-3}}, & 0 \leqslant r < \SI{9}{\micro \metre},  \\
         \left(1 + 40 \left(\frac{r}{\SI{2}{\micro \metre}}\right)^7-140\left(\frac{r}{\SI{2}{\micro \metre}}\right)^6+168\left(\frac{r}{\SI{2}{\micro \metre}}\right)^5-70\left(\frac{r}{\SI{2}{\micro \metre}}\right)^4\right)\times 10^{17} {\rm cm^{-3}}, & \SI{9}{\micro \metre} \leqslant r \leqslant \SI{11}{\micro \metre}, \\
         -10^{17} {\rm cm^{-3}}, & r > \SI{11}{\micro \metre}
    \end{array}
     \\
    \right.
\end{equation}
where $\textbf{x} = (x,y)^{\intercal}$ and the radial parameter $r = \sqrt{x^2 + y^2}$.
$\hat{N}(r)$ is plotted in Fig. \ref{fig:Exp_3_2D_case} for illustration. 
The cathode and anode electrodes are Dirichlet boundaries for $\psi,n,p$ and others are Neumann boundaries.

\begin{figure}[!h]
    \begin{subfigure}[b]{0.5\textwidth}
        \centering
        \includegraphics[width=0.6\textwidth]{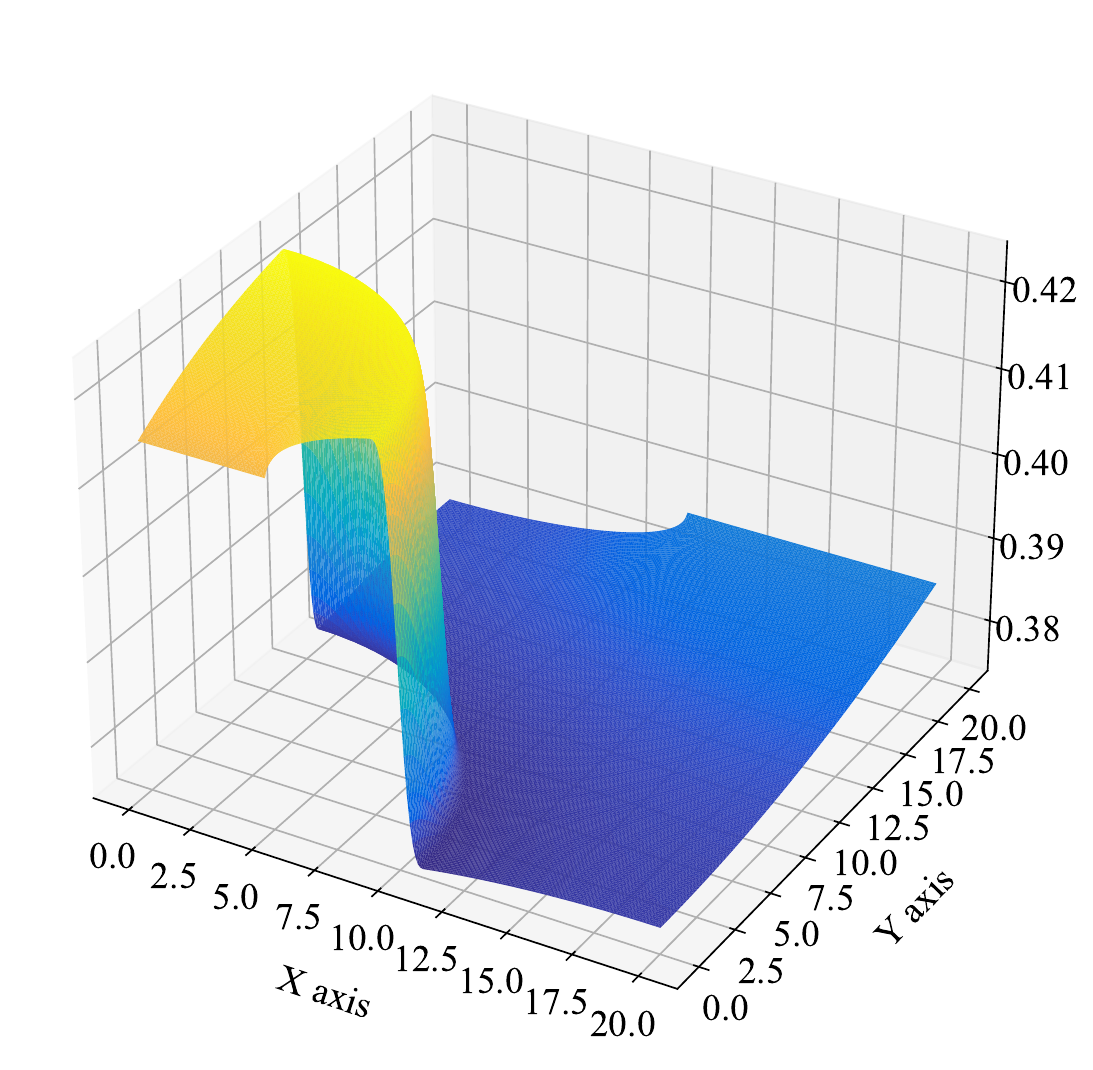}
        \label{fig:Exp_3_reference_solution_psi}
    \end{subfigure}
    \begin{subfigure}[b]{0.5\textwidth}
        \centering
        \includegraphics[width=0.6\textwidth]{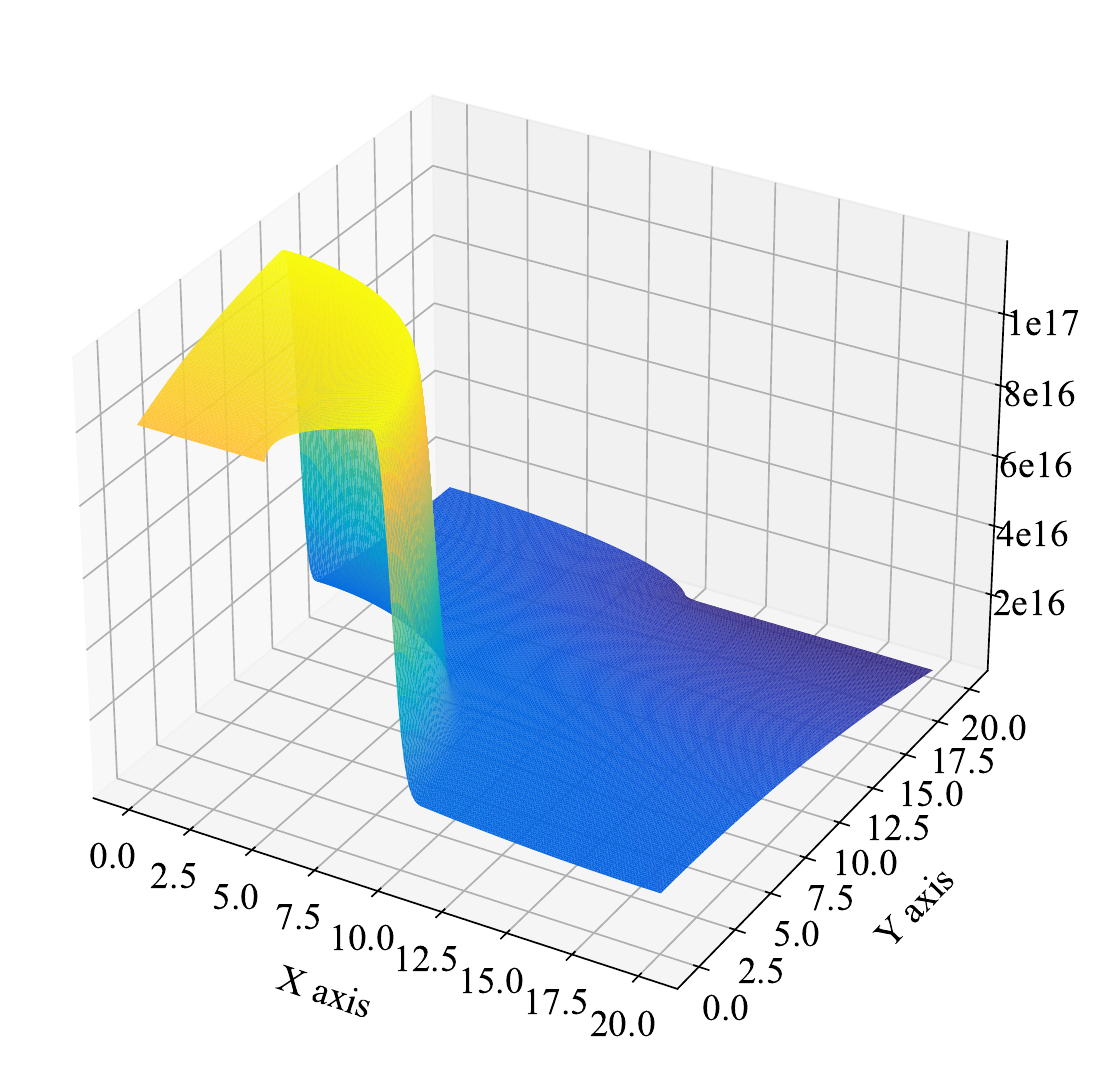}
        \label{fig:Exp_3_reference_solution_n}
    \end{subfigure}
    \caption{Solutions of the general 2D test case in subsection \ref{sec:num_ha_indicator}: $\psi$ (left) and $n$ (right).}
    \label{fig:Exp_3_reference_solution}
\end{figure}

This example simulates the conducting steady state of the diode, with the anode electrode forward biased at $0.8 \rm V$ and the cathode electrode grounded at $0 \rm V$. As shown in Fig. \ref{fig:Exp_3_reference_solution}, both $\psi$ and $n$ exhibit sharp layers near the junction, where the numerical stability of conventional HDG schemes is hindered.

As demonstrated in the previous examples, the HA technique shall be used in the junction regions to ensure the stability of HDG schemes. Therefore, effective indicators should identify the cells in these regions where the conventional HDG scheme may encounter convergence issues. To this aim, we propose three indicators as stated in subsection \ref{sec:ha-adaptivity}, namely $\|\nabla \psi \|_{L^2(K)},\|\nabla N \|_{L^2(K)}$ and $\|\nabla n \|_{L^2(K)}$ for each cell $K$. They are computed in this test case to validate their effectiveness in identifying cells near junctions. 

\begin{figure}[!ht]
    \centering
    \begin{subfigure}[b]{0.33\textwidth}
        \centering
        \includegraphics[width=1.0\textwidth]{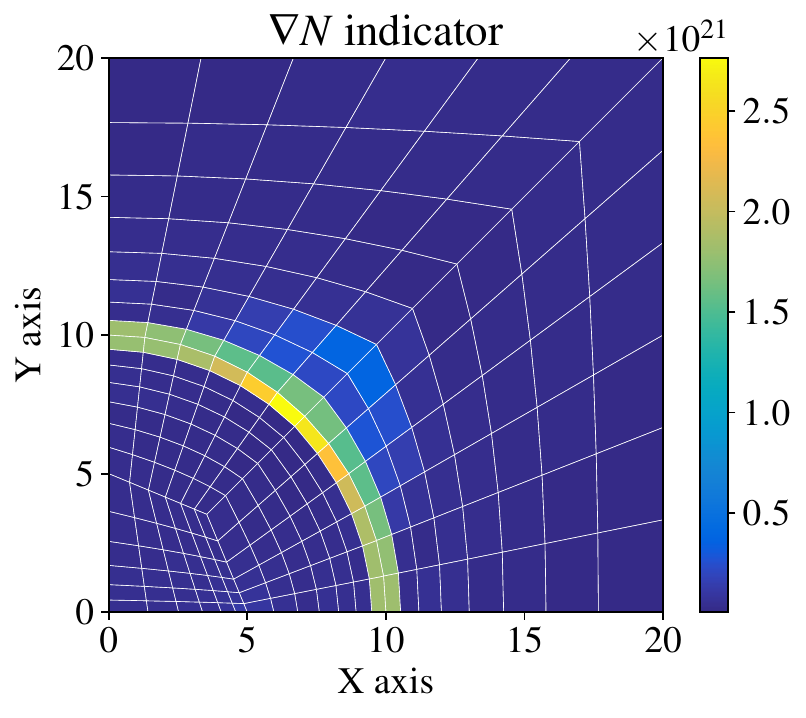}
        \label{fig:Exp_3_grad_doping_indicator}
    \end{subfigure}
    \begin{subfigure}[b]{0.33\textwidth}
        \centering
        \includegraphics[width=1.0\textwidth]{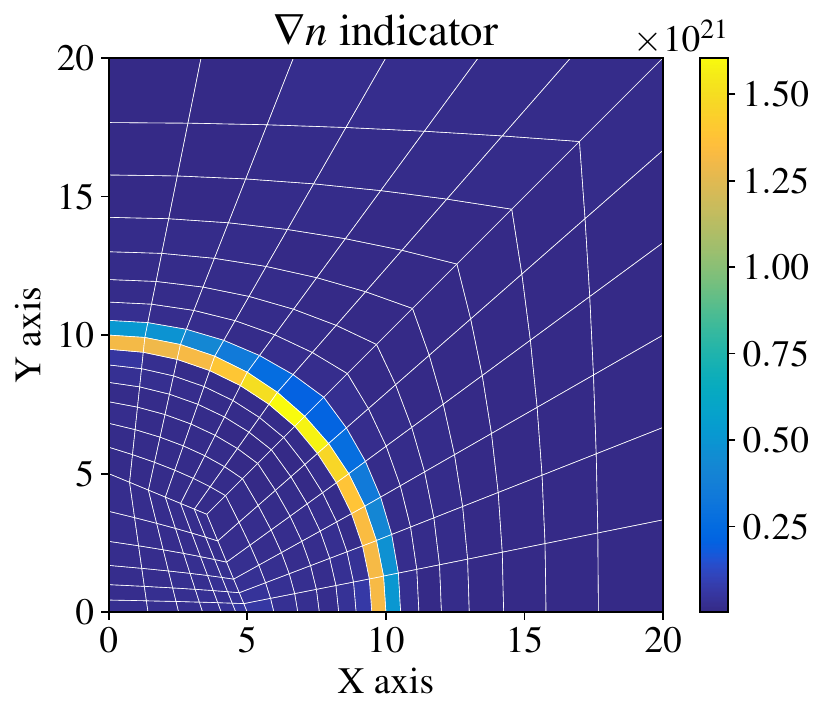}
        \label{fig:Exp_3_grad_n}
    \end{subfigure}
    \begin{subfigure}[b]{0.33\textwidth}
        \centering
        \includegraphics[width=1.0\textwidth]{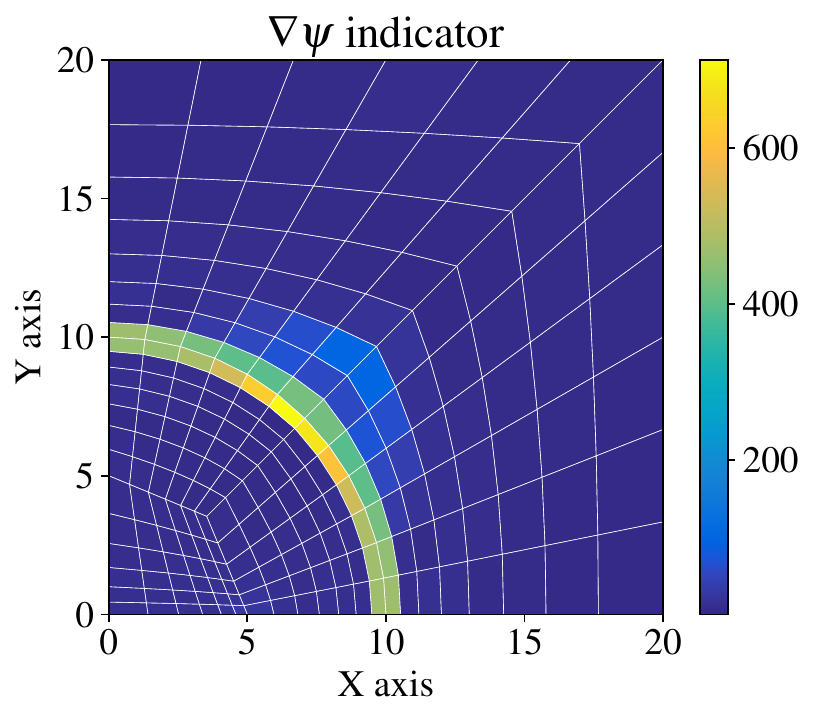}
        \label{fig:Exp_3_grad_psi}
    \end{subfigure}
    \caption{Indicators for HA-adaptive: $\|\nabla N\|_{L^2(K)}$ (left), $\|\nabla n\|_{L^2(K)}$ (middle) and $\|\nabla \psi\|_{L^2(K)}$ (right).}
    \label{fig:Exp_3_different_HA_indicators}
\end{figure}

The HDG-HA scheme is applied to all cells to obtain a convergent result. Three indicators are then calculated for each cell, producing cell-wise constant values, as shown in Fig. \ref{fig:Exp_3_different_HA_indicators}. It is clear that all three proposed indicators can accurately identify the junction region. Nevertheless, we prefer using $\| \nabla \psi \|_{L^2(K)}$-based indicator in practical implementation because the electric field $\mathbf{E}=-\nabla\psi$ is the intrinsic source of the convection-dominant instability. 

\begin{figure}[!ht]
    \centering
    \centering
    \includegraphics[width=0.3\textwidth]{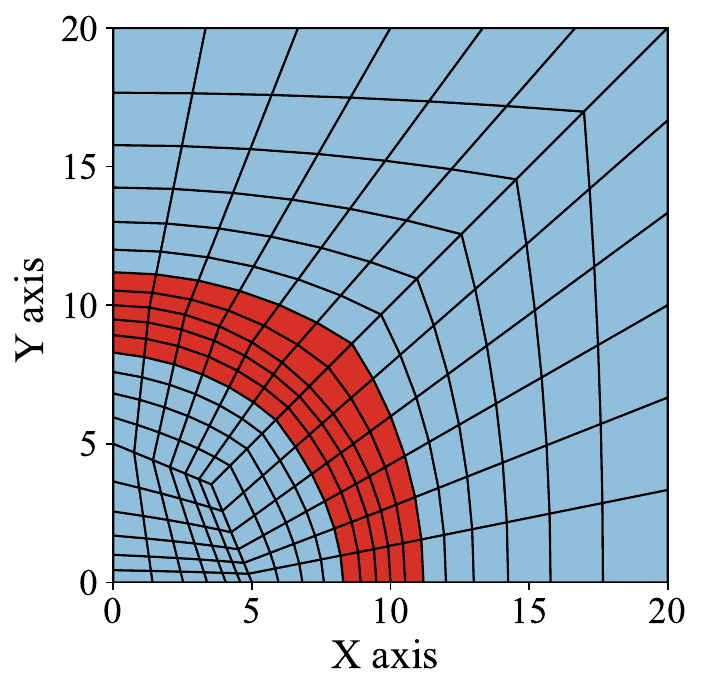}
    \caption{Conventional HDG cells (in blue) and HDG-HA cells (in red).}
    \label{fig:Exp_3_HA_indicator}
\end{figure}

In the remainder of this paper, to apply the HA-adaptive procedure, we designate cells $K_i$ as HA cells where the indicator $\| \nabla \psi \|_{L^2(K_i)}>20\%\max\limits_{K_j\in\mathcal{T}_h}\| \nabla \psi \|_{L^2(K_j)}$. The HDG-HA scheme is applied to the HA cells, while the conventional HDG schemes are applied to the other cells. This selection process is illustrated in Fig. \ref{fig:Exp_3_HA_indicator}, showing that the HA cells correspond well with the junction region in Fig. \ref{fig:Exp_3_2D_case}. 

\subsection{$hp$-adaptivity}

In addition to HA-adaptivity, which switches between HA cells and conventional HDG cells, it should be further remarked that $hp$-adaptivity can be seamlessly integrated into our HDG-HA framework.

\begin{figure}[!ht]
    \begin{subfigure}[b]{0.38\textwidth}
        \centering
        \includegraphics[width=1\textwidth]{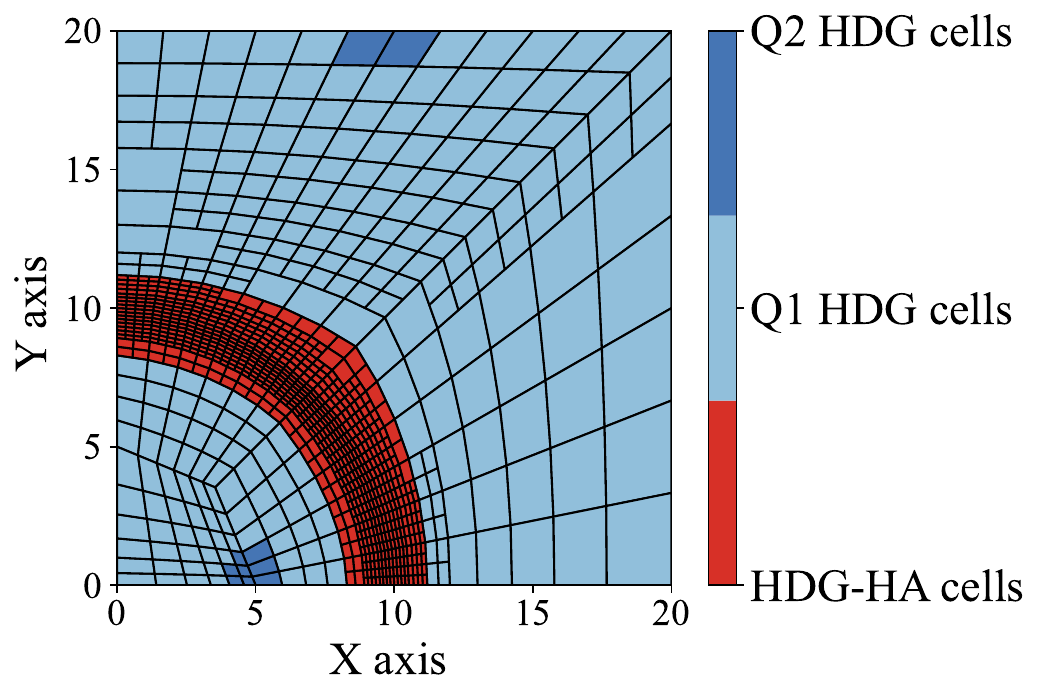}
        \label{fig:Exp_3_HA_hp_1}
    \end{subfigure}
    \begin{subfigure}[b]{0.3\textwidth}
        \centering
        \includegraphics[width=1\textwidth]{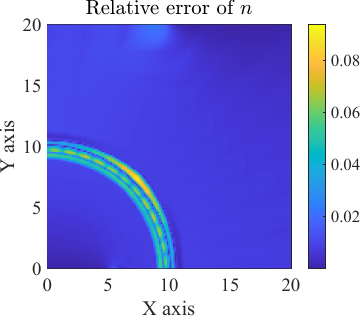}
        \label{fig:Exp_3_HA_hp_2}
    \end{subfigure}
    \begin{subfigure}[b]{0.3\textwidth}
        \centering
        \includegraphics[width=1\textwidth]{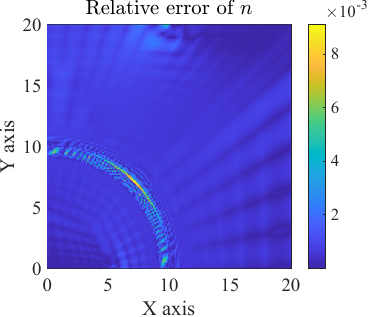}
        \label{fig:Exp_3_HA_hp_3}
    \end{subfigure}
    \caption{Results of $hp$-adaptive on the 2D diode test case: the mesh of HA- and $hp$-adaptive (left),  the error of applying only HA-adaptivity (middle) and the error of applying both HA- and $hp$-adaptivity (right).}
    \label{fig:Exp_3_HA_hp_adaptive}
\end{figure}

With the stability greatly enhanced by the HA-adaptive process (Fig. \ref{fig:Exp_3_HA_indicator}), the numerical accuracy of the HDG-HA scheme can be further improved by $h$-adaptivity illustrated in Fig. \ref{fig:Exp_3_HA_hp_adaptive}.
Based on the HA-adaptivity shown in Fig. \ref{fig:Exp_3_HA_indicator}, we further conduct two additional levels of $h$-refinement for the cells in the junction region and near the top boundary, as clearly shown in Fig. \ref{fig:Exp_3_HA_hp_adaptive}.
The $p$-refinement is also conducted on regions where error is still big and the solution is smooth. 
After the $hp$-refinement, the maximum relative error of $n$ is significantly reduced by 1 order of magnitude while the number of unknowns only increase by 4 times, compared to the single HA-adaptive procedure in Fig. \ref{fig:Exp_3_HA_indicator}. 
\rtwo{
Additionally, even after applying two additional levels of $h$-refinement, the maximum error still occurs in the HA region. This suggests that, for a fixed number of cells, the HA region should be refined further to maintain relatively high accuracy over the entire domain. As an initial exploration of the HDG-HA scheme, we demonstrate the implementation of $hp$-adaptivity in this simple example. Its significant potential for achieving high accuracy with fewer degrees of freedom should be further investigated in future work.}

\section{Real examples from industry}
\label{sec:RealWorldExamples}

\subsection{A 2D P-type-intrinsic-N-type (PiN) power diode}

\begin{figure}[!ht]
    \centering
    \begin{subfigure}[b]{0.43\textwidth}
        \centering
        \includegraphics[width=0.64\textwidth]{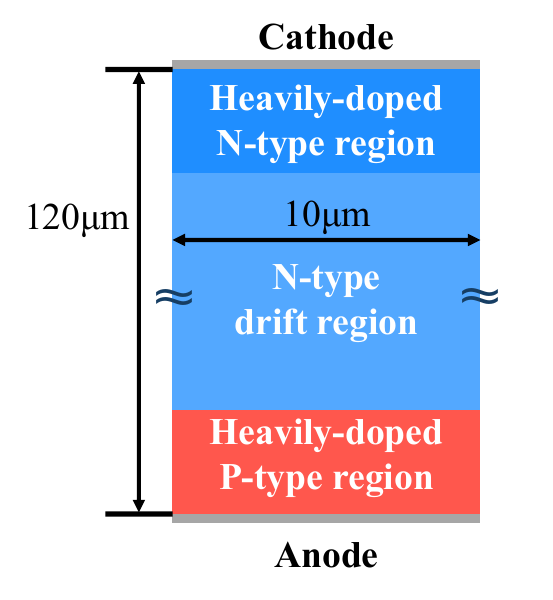}
        \label{fig:Exp_3_p_initial_mesh}
    \end{subfigure}
    \begin{subfigure}[b]{0.4\textwidth}
        \centering
        \includegraphics[width=1\textwidth]{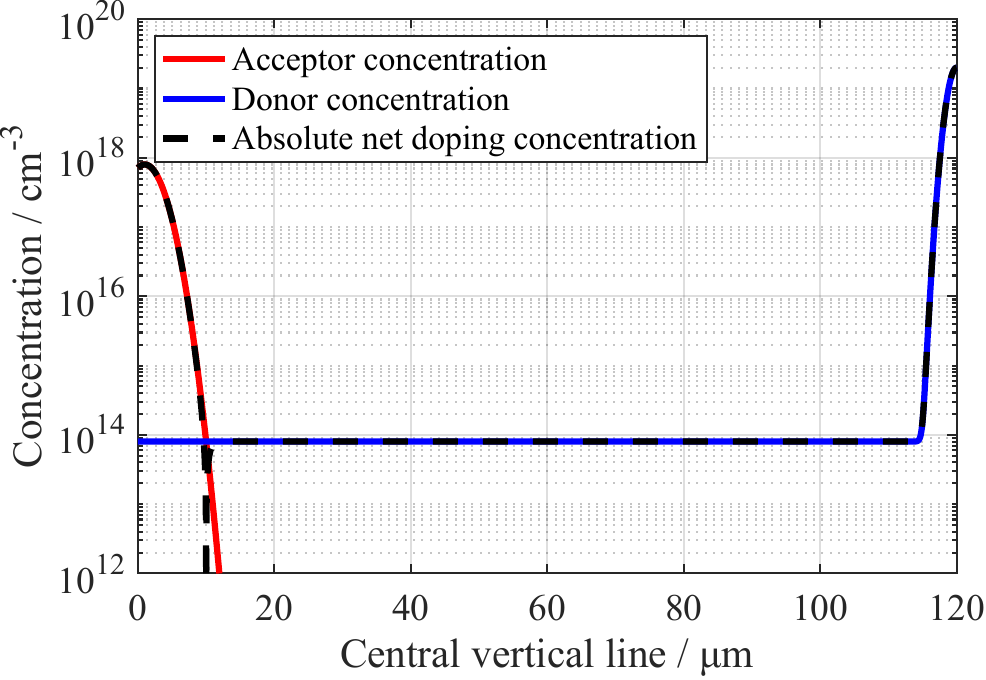}
        \label{fig:Exp_3_p_indicator}
    \end{subfigure}
\caption{The structure and doping of the power PiN diode: structure illustration (left) and doping profile (right).}
\label{fig:Exp_4_general_diode_structure}
\end{figure}

We first employ the proposed HDG-HA scheme to simulate a realistic power diode to validate its stability and accuracy. A PiN diode \cite{Baliga2010Book}, characterized by an intrinsic (or lightly-doped) region sandwiched between a heavily-doped P-type and a heavily-doped N-type region, is selected as the test case. This device is typically fabricated by diffusing dopants from opposite sides into a lightly doped bulk substrate. Fig. \ref{fig:Exp_4_general_diode_structure} illustrates the structure of the power diode and its doping profile. The size of the diode is $\SI{120}{\micro \metre}\times\SI{10}{\micro \metre}$, and the net doping concentration changes sharply from $-10^{16}\rm cm^{-3}$ to $8\times 10^{13} \rm cm^{-3}$ within only $\SI{4}{\micro \metre}$, forming a sharp PN-junction, as shown in Fig. \ref{fig:Exp_4_general_diode_structure}. The two electrodes on the top and bottom sides of the device, i.e., cathode and anode, are two Dirichlet boundaries for all solution variables $\psi,n,p$. The left and right sides are Neumann boundaries for $\psi,n,p$ with no fluxes. The independent scaling bases are $x^*=\SI{120}{\micro \metre},V^*=0.026{\rm V},N^*=2\times 10^{19} \rm cm^{-3},\mu^*=1417 \rm cm^2/(V s)$, resulting in $\lambda^2 \approx 10^{-10}$.

\begin{figure}[!ht]
\centering
\begin{subfigure}[b]{0.32\textwidth}
    \includegraphics[width=1\textwidth]{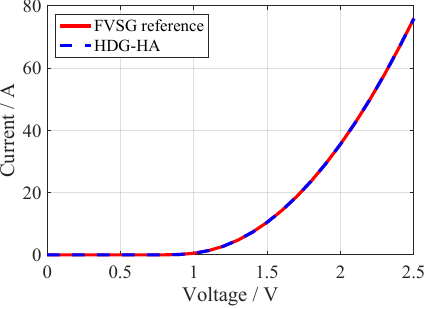}
\end{subfigure}
\begin{subfigure}[b]{0.33\textwidth}
    \includegraphics[width=1\textwidth]{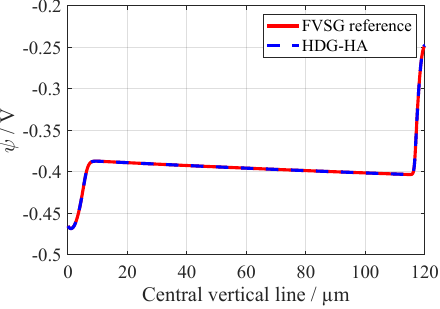}
\end{subfigure}
\begin{subfigure}[b]{0.33\textwidth}
    \includegraphics[width=1\textwidth]{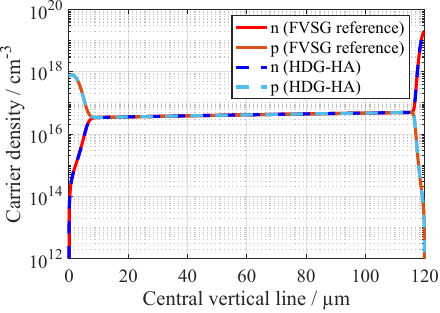}
\end{subfigure}

\caption{\rone{Current-voltage characteristic and the $0.8 \rm V$ solutions of the PiN diode: current-voltage characteristic (left), potential (middle) and carrier concentrations (right). 8.64$\times10^6$ cells are used in FVSG reference, while HDG-HA scheme uses 1.92$\times10^4$ cells.}}
\label{fig:Exp_4_general_diode_IV_curve}
\end{figure}

The current-voltage characteristic of this diode is simulated by the proposed HDG-HA scheme. The thickness of the diode is $10^7\si{\micro \metre}$, and we assume all physical quantities are homogeneous in the third dimension. Therefore the 2D areas of the cathode and anode electrodes are both $1 \rm cm^2$. As in Fig. \ref{fig:Exp_4_general_diode_IV_curve}, the result of HDG-HA scheme closely matches the reference FVSG result obtained from an extremely fine mesh. In Fig. \ref{fig:Exp_4_general_diode_IV_curve}, the unscaled $V_{\rm applied}$ for cathode boundary stays $0 \rm V$ while the unscaled $V_{\rm applied}$ for anode boundary increases along the horizontal axis. \rone{The solutions along the central vertical line also match the reference FVSG results without any oscillations.}

\subsection{A 2D power bipolar junction transistor (BJT)}

\begin{figure}[!ht]
    \begin{subfigure}[b]{0.33\textwidth}
        \centering
        \includegraphics[width=0.72\textwidth]{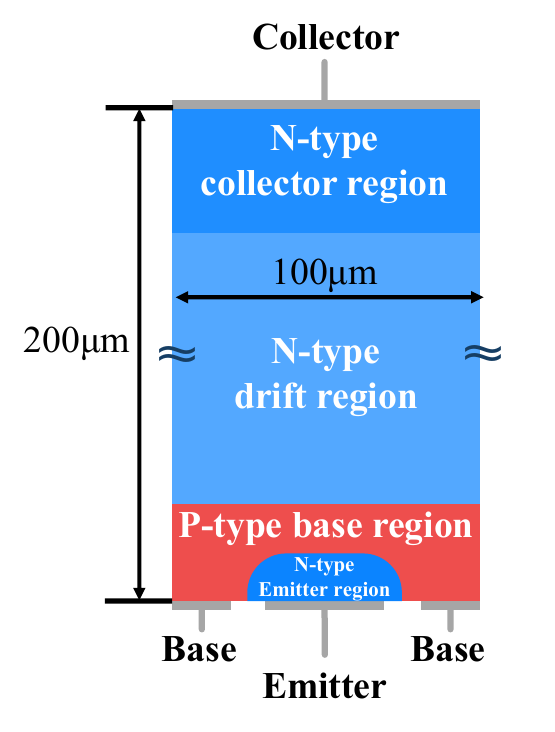}
    \end{subfigure}
    \begin{subfigure}[b]{0.33\textwidth}
        \centering
        \includegraphics[width=1.0\textwidth]{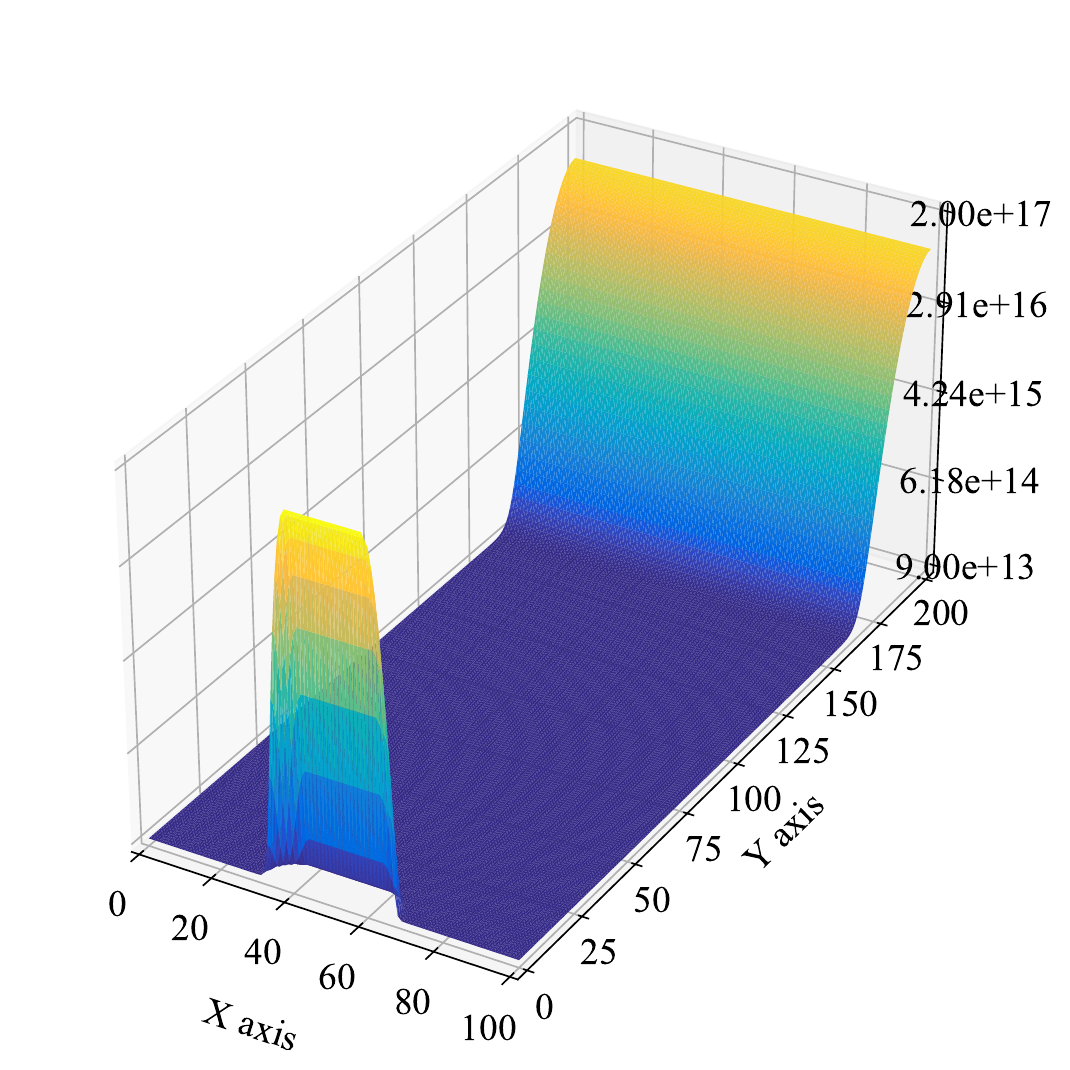}
    \end{subfigure}
    \begin{subfigure}[b]{0.33\textwidth}
        \centering
        \includegraphics[width=1.0\textwidth]{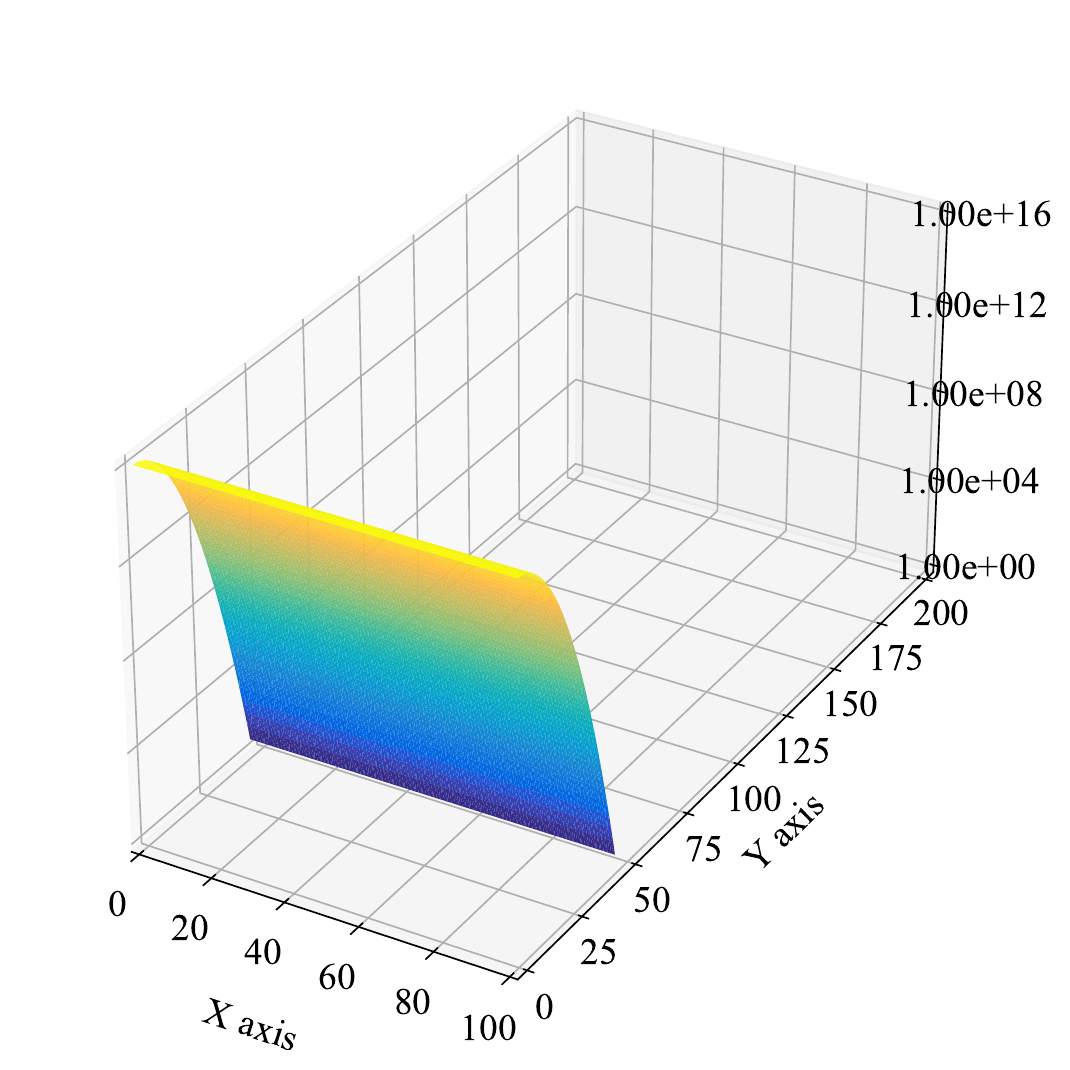}
    \end{subfigure}
\caption{The doping and structure of a power BJT: structure illustration (left), donor concentration (middle) and acceptor concentration (right).}
\label{fig:Exp_5_BJT_structure}
\end{figure}

The bipolar junction transistors (BJTs)  \cite{Baliga2010Book} are a family of power semiconductor devices widely used in power electronic systems with the advantage of high voltage-blocking and current-surging ability. It has three metal electrodes and two PN-junctions, and its manufacturing involves multiple processes. A real-world power BJT in Fig. \ref{fig:Exp_5_BJT_structure} is taken as an example. The simulated 2D cross section of the BJT is $\SI{200}{\micro \metre} \times \SI{100}{\micro \metre}$ in size with several regions:

\begin{enumerate}
    \item Drift region: The wide drift region comes from the original silicon epitaxial wafer, which has a uniform Phosphorus concentration of $\rm 9 \times 10^{13} cm^{-3}$;

    \item Base region: Manufactured by ion implantation and diffusion process with Boron. The concentration of dopant is a Gaussian distribution with $\rm 10^{16} cm^{-3}$ maximum;

    \item Collector region: Manufactured by ion implantation and diffusion process with Phosphorus. The concentration of dopant is a Gaussian distribution with $\rm  10^{17} cm^{-3}$ maximum;

    \item Emitter region: Manufactured by ion implantation and diffusion process with Phosphorus.  The concentration of dopant is a Gaussian distribution with $\rm 2 \times 10^{17} cm^{-3}$ maximum.
\end{enumerate}

For this BJT device, the Base, Emitter, and Collector electrode contacts are Dirichlet boundaries, whereas the remaining boundaries are the Neumann boundaries. More parameters of this BJT device can be found in the Appendix. The independent scaling bases are chosen as $x^*=\SI{200}{\micro \metre},V^*=0.026V,N^*=2\times 10^{17} \rm cm^{-3},\mu^*=1417 \rm cm^2/(V s)$, resulting in $\lambda^2 \approx 2\times 10^{-9}$. 

In this experiment, the HDG-HA scheme is employed to simulate and analyse the power BJT under voltage blocking working conditions. When high voltage is applied across the Collector-Emitter electrodes, huge electric field will be built up in the drift region with a punch-through effect. The punch-through is a typical mechanism of large size power BJTs that improves voltage blocking ability without sacrificing device thickness. In a one dimensional simplified model, the total blocked voltage across Collector-Emitter is the integral of electric field $V_{blocked} = \int_{Collector}^{Emitter} |\mathbf{E}| dx$, which means that maximizing electric filed $|\mathbf{E}|$ along the integration line could increase blocking voltage with the same device thickness (distance from Collector to Emitter).

\begin{figure}[!ht]
    \begin{subfigure}[b]{0.5\textwidth}
        \centering
        \includegraphics[width=1\textwidth]{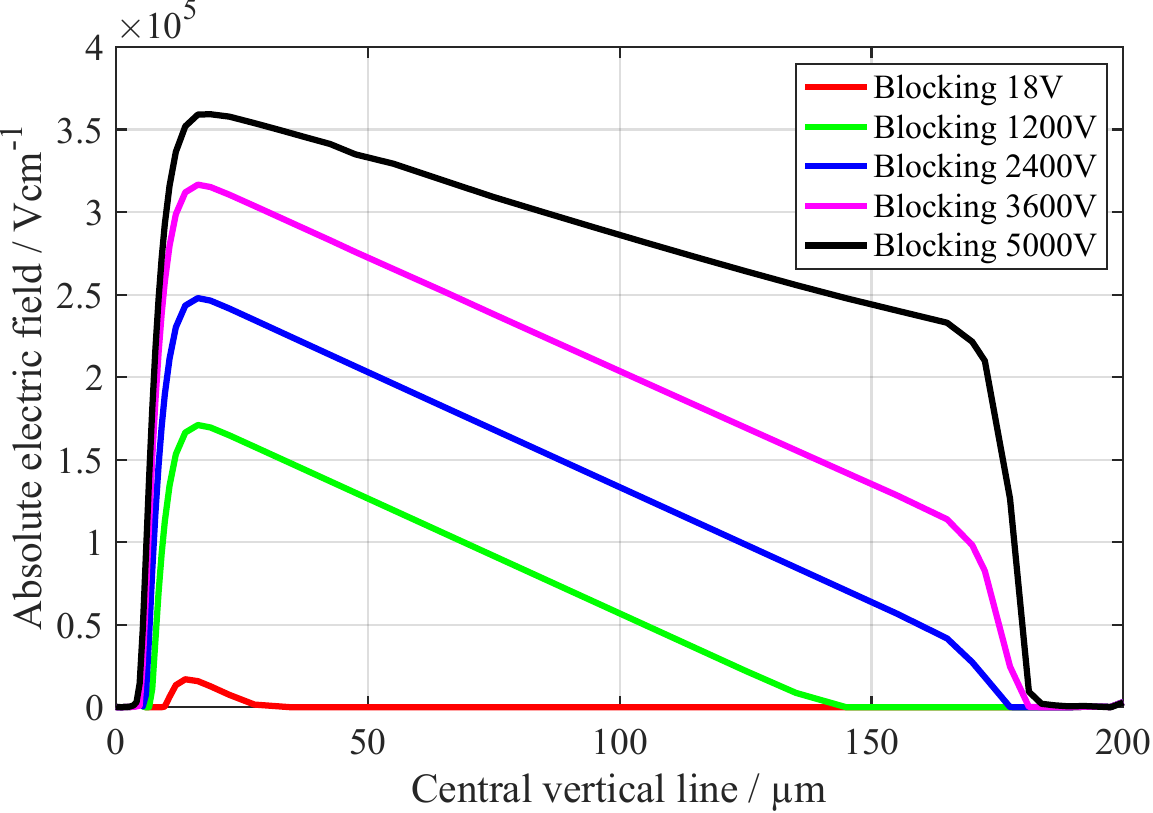}
    \end{subfigure}
    \begin{subfigure}[b]{0.5\textwidth}
        \centering
        \includegraphics[width=0.76\textwidth]{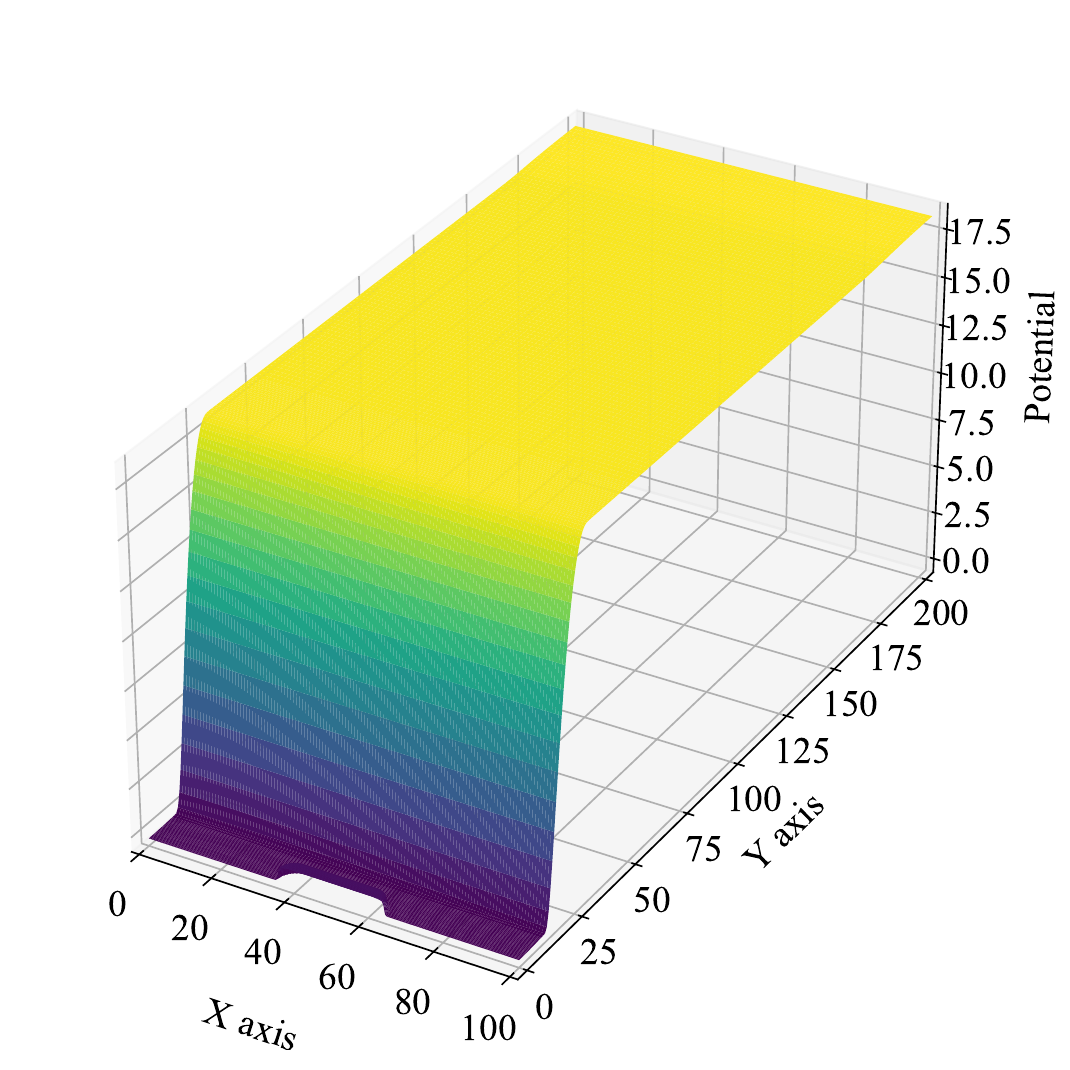}
    \end{subfigure}
\caption{Illustration the punch-through electric field in a BJT device (left) and electrostatic potential $\psi$ at $18\rm V$ (right).}
\label{fig:Exp_5_BJT_punch_through}
\end{figure}

We use the proposed HDG-HA scheme to analyze the voltage blocking capability of the BJT device. Different voltages across the Collector and the Emitter electrodes are simulated, the absolute electric field on the line across the center of Collector and Emitter electrodes is shown in Fig. \ref{fig:Exp_5_BJT_punch_through}. In this experiment, the unscaled $V_{\rm applied}$ for Base and Emitter electrodes are $0\rm V$ whereas the unscaled $V_{\rm applied}$ for the Collector electrode is shown in Fig. \ref{fig:Exp_5_BJT_punch_through}. At lower blocking voltages, the distribution of electric field strength presents a triangular shape, confined to the sides with high and low doping levels, respectively. However, at higher blocking voltages, the electric field distribution spans the entire low-doped region, exhibiting an approximately trapezoidal shape. Such a distribution of the electric field enables this power device to block higher voltages without increasing thickness. The HDG-HA scheme successfully depicts this punch-through mechanism.

\section{Conclusion}
We proposed a hybridizable discontinuous Galerkin (HDG) scheme with harmonic averaging (HA) technique to simulate the real-world semiconductor devices. The proposed HDG-HA scheme combines the robustness of finite volume Scharfetter-Gummel method with the high-order accuracy (in the areas with smooth solutions) and $hp$-flexibility offered by the locally conservative HDG scheme. The coupled Poisson equation and two drift-diffusion equations are simultaneously solved by the Newton method. Indicators based on the gradient of net doping and solution variables are proposed to act as a switch sign of using HA technique. Numerical results showed that the proposed scheme experiences no oscillations or convergence issues even on heavily doped sharp PN-junctions. Devices with circular junctions and realistic doping profiles have been simulated in two dimensions, qualifying this scheme for the practical simulation of real-world semiconductor devices. 

\section*{Acknowledgments}
This work is partly supported by Independent Innovation Research Program of Department of Electrical Engineering in Tsinghua University (C.~Zhuang), the NSFC grants 12171041 and 11771036 (Y.~Cai) and the Ministry of Education of Singapore under its AcRF Tier 2 funding MOE-T2EP20122-0002 (A-8000962-00-00) (B.~Lin and W.~Bao).

\bibliographystyle{model1-num-names}
\bibliography{refs}

\section{Appendix}
\label{sec:Appendix}

In this paper, the physical quantities are chosen as the average constant values of silicon material, as listed in Table \ref{tab:Appendix1}.

\begin{table}[!ht]
\caption{\label{tab:Appendix1}Basic physical quantities}
\centering
\begin{tabular}{|l|l|l|}
\hline
Physical quantity  & Symbol & Numerical value \\ \hline
Elementary charge  &  $q$   & $1.60217663\times 10^{-19} \ \rm C$ \\ \hline
Boltzmann constant &  $k$   & $1.38064852\times 10^{-23} \ {\rm J/K}$ \\ \hline
Thermal voltage &  $V_T$   &  $0.02585199 \ \rm V$ \\ \hline
Silicon dielectric constant & $\varepsilon$ & $1.03593997\times10^{-12} \ \rm C/(V \cdot cm)$  \\ \hline
Effective intrinsic carrier concentration &  $n_{\rm ie}$   & $1.08738184\times 10^{10}  \ \rm cm^{-3}$  \\ \hline
Electron mobility &  $\mu_n$  & $1417.0 \ {\rm cm^2/(V \cdot s)}$  \\ \hline
Hole mobility &  $\mu_p$      & $470.5 \ {\rm cm^2/(V \cdot s)}$ \\ \hline
Electron diffusion coefficient &  $D_n$  & $36.63227105 \ \rm cm^2/s $  \\ \hline
Hole diffusion coefficient &  $D_p$      & $12.16336170 \ \rm cm^2/s $ \\ \hline
Electron lifetime &  $\tau_n$      & $1\times 10^{-3} \ \rm s$   \\ \hline
Hole lifetime &  $\tau_p$      & $3\times 10^{-4} \ \rm s$ \\ \hline
Electron recombination coeﬃcients &  $C_n$  & $6.59841820\times 10^{-31} \ \rm cm^6/s$  \\ \hline
Hole recombination coeﬃcients &  $C_p$  & $4.15058741\times 10^{-31} \ \rm cm^6/s$ \\ \hline
\end{tabular}
\end{table}

The basic units of scaling in section \ref{subsec:Nondimensionless} are listed in Table \ref{tab:Appendix2}

\begin{table}[!ht]
\caption{\label{tab:Appendix2}Scaling basic units}
\centering
\begin{tabular}{|l|l|l|l|}
\hline
Used for physical quantity  & Symbol & Choice  & Typical value \\ \hline
Carrier concentrations & $N^*$ & $\max(N(x))$  & $10^{19} \rm cm^{-3}$ \\ \hline
Voltages & $V^*$ & $V_T$  & $0.02585199 \ \rm V$ \\ \hline
Space dimension & $x^*$ & $\max(\Omega)$  & $ \SI{1000}{\micro \metre}$ \\ \hline
Carrier diffusion coefficient & $D^*$ & $\max(D_n,D_p)$  & $36.63227105 \ \rm cm^2/s $\\ \hline
Carrier mobility & $\mu^*$ & $D^*/V^*$  & $1417.0 \ {\rm cm^2/(V \cdot s)}$\\ \hline
Time & $t^*$ & $(x^*)^2/D^*$  & $ 2.72983348 \times 10^{-4} \ \rm s$ \\ \hline
Current densities & $J^*$ & $q D^*N^*/x^*$  & $5.86919316 \times 10^{2} \ \rm A/cm^2$ \\ \hline
Recombination rates & $R^*$ & $D^*N^*/(x^*)^2$  & $3.66322711 \times 10^{22} \ \rm cm^{-3}/s$ \\ \hline
Recombination coeﬃcients & $C^*$ & $1/((N^*)^2t^*)$  & $3.66322711 \times 10^{-35} \ \rm cm^{6}/s$ \\ \hline
\end{tabular}
\end{table}

The doping profile of the PiN diode test case is made from three process steps:

\begin{enumerate}
    \item Uniform Phosphorus doping of $8\times 10^{13} \rm cm^{-3}$,

    \item Gaussian doping of Boron. Boron concentration is $8 \exp\left(-\left(\frac{y-\SI{1}{\micro\metre}}{\SI{2.97}{\micro\metre}}\right)^2\right) \times 10^{17} \rm cm^{-3}$ where $y$ is the unscaled vertical coordinate,

    \item Gaussian doping of Phosphorus. Phosphorus concentration is $2 \exp\left(-\left(\frac{\SI{120}{\micro \metre} - y}{\SI{1.42}{\micro\metre}}\right)^2\right) \times 10^{19} \rm cm^{-3}$ where $y$ is the unscaled vertical coordinate.

\end{enumerate}

The doping profile of the power bipolar junction transistor (BJT) test case is made from 4 process steps:

\begin{enumerate}
    \item Uniform Phosphorus doping of $9\times 10^{13} \rm cm^{-3}$. This makes up the drift region;

    \item Gaussian doping of Boron. Boron concentration is $ \exp\left(-\left(\frac{y}{\SI{7.5}{\micro \metre}}\right)^2\right) \times 10^{16} \rm cm^{-3}$ where $y$ is the unscaled vertical coordinate. This makes up the P-type base region;

    \item Gaussian doping of Phosphorus. Phosphorus concentration is $ \exp\left(-\left(\frac{\SI{200}{\micro \metre} - y}{\SI{12}{\micro\metre}}\right)^2\right) \times 10^{17} \rm cm^{-3}$ where $y$ is the unscaled vertical coordinate. This makes up the N-type Collector region;
    
    \item Gaussian doping of Phosphorus. Phosphorus concentration is $ 2\exp\left(-\left(\frac{d}{\SI{2.5}{\micro\metre}}\right)^2\right) \times 10^{17} \rm cm^{-3}$ where $d$ is the unscaled distance to the closest point on the doping baseline. The doping baseline is line segment $(\SI{40}{\micro\metre}, \SI{0}{\micro \metre})$ to $(\SI{60}{\micro\metre}, \SI{0}{\micro \metre})$. This makes up the N-type Emitter region with a rounded boundary.

\end{enumerate}

And the positions of the three electrodes of the BJT are:

\begin{enumerate}
    \item The Base electrode is made of two parts connected by external circuit, which means that the two parts share the same Dirichlet boundary condition for $\psi, n, p$. The first part is from $(\SI{0}{\micro \metre},\SI{0}{\micro \metre})$ to $(\SI{20}{\micro \metre},\SI{0}{\micro \metre})$ and the second part is from $(\SI{80}{\micro \metre},\SI{0}{\micro \metre})$ to $(\SI{100}{\micro \metre},\SI{0}{\micro \metre})$;
    
    \item The Emitter electrode is from $(\SI{40}{\micro \metre},\SI{0}{\micro \metre})$ to $(\SI{60}{\micro \metre},\SI{0}{\micro \metre})$;

    \item The Collector electrode is from $(\SI{0}{\micro \metre},\SI{200}{\micro \metre})$ to $(\SI{100}{\micro \metre},\SI{200}{\micro \metre})$.
\end{enumerate}

\end{document}